\documentclass[11pt]{amsart}
\usepackage[utf8]{inputenc}

\usepackage{amsmath,amsthm,amssymb}
\usepackage{amsxtra} 
\usepackage{amsfonts}

\usepackage{txfonts}
\usepackage{breakurl}
\usepackage{bbm}
\usepackage{mathrsfs}
\usepackage{xstring}
\usepackage{diagbox}
\usepackage[mathcal]{euscript}

\usepackage{enumitem}

\usepackage[a4paper]{geometry}

\usepackage{float}

\usepackage{graphicx}
\usepackage{upref}
\usepackage{color}
\usepackage[active]{srcltx}
\usepackage[colorlinks=true,linkcolor=blue,citecolor=green,urlcolor=red]{hyperref}
\usepackage[usenames,dvipsnames]{pstricks} 
\usepackage{epsfig} 
\usepackage{subcaption}
\usepackage{tikz}

\usepackage{color}
\usepackage{tcolorbox}

\newtheorem{thm}{Theorem}[section]

\newtheorem{Cor}[thm]{Corollary}

\theoremstyle{definition}
\newtheorem{Def}[thm]{Definition}

\newtheorem{remark}[thm]{Remark}

\newtheorem{Ej}[thm]{Example}

\numberwithin{equation}{section}

\DeclareMathAlphabet{\pazocal}{OMS}{zplm}{m}{n}

\newcommand{\Z}{\mathbb{Z}}

\newcommand{\K}{\mathbb{K}}

\newcommand{\R}{\mathbb{R}}
\newcommand{\Sy}{\mathbb{S}}

\newcommand{\Rep}{\operatorname{Rep}}
\newcommand{\Traza}{\operatorname{tr}}

\newcommand{\ep}{\varepsilon}

\DeclareMathOperator{\Hom}{Hom}

\newcommand{\g}{\mathfrak{g}}

\newcommand{\h}{\mathfrak{h}}

\newcommand{\Sl}{\mathfrak{sl}}
\newcommand{\Gl}{\mathfrak{gl}}

\newcommand{\SL}{\mathrm{SL}}

\newcommand{\A}{\mathfrak{A}}
\newcommand{\B}{\mathfrak{B}}

\newcommand{\bas}{\mathcal{B}}

\newcommand{\bwedge}{\mbox{\normalsize $\bigwedge$}}

\newcommand{\hdot}{{\raise1.2pt\hbox to0.3em{\Large $\dot$}}}

\begin{document}

\title{$G$-tables and the Poisson structure of the even
cohomology of cotangent bundle of the Heisenberg Lie group}

\author{Leandro Cagliero}
\address{FaMAF-CIEM (CONICET), Universidad Nacional de C\'ordoba,
Medina Allende s/n, Ciudad Universitaria, 5000 C\'ordoba, Rep\'ublica Argentina.}
\email{cagliero@famaf.unc.edu.ar}

\author{Gonzalo Gutierrez}
\address{CIEM-CONICET, FCEyN-UBA, Intendente Güiraldes 2160 - Ciudad Universitaria, CABA, Rep\'ublica Argentina.}
\email{gonzalounsafcex@gmail.com}

\thanks{This research was partially supported by a CONICET grant
PIP 11220210100597CO, and SeCyT-UNC grant 33620230100117CB}

\subjclass[2020]{17B05, 17B63, 22E46, 20F29}

\keywords{non-semisimple Lie algebras, uniserial representations, tensor product, socle, radical, intertwining operators}

\begin{abstract}
    In the first part of the paper, 
     we define the concept of a $G$-table of a 
     $G$-(co)algebra and we compute the $G$-table of some $G$-(co)algebras 
     (here a $G$-algebra is an algebra on which $G$ acts, semisimply, by algebra automorphisms). 
     The $G$-table of a 
     $G$-(co)algebra $\A$ is a set of scalars that provides very precise and concise information about both the algebra structure and the $G$-module structure of $\A$. In particular, the ordinary multiplication table of $\A$ can be derived from the $G$-table of $\A$.
     From the 
    $G$-table of a $G$-algebra $\A$ 
    we define a plain algebra $P(\A)$ associated to it and we present some basic functoriality results about $P$.

Obtaining the $G$-table of a given $G$-algebra $\A$ requires a considerable amount of work but, the result, is a very powerful tool as 
shown in the second part of the paper. 
    Here we compute the 
    $\SL(2,\K)$-tables of the Poisson algebra structure of
    the even-degree part of the cohomology associated to the cotangent bundle of the 3-dimensional Heisenberg Lie group with Lie algebra $\h$, that is 
    $H_E^{\bullet,\bullet}(\h)=H_E^{\bullet}(\h,\bigwedge^{\bullet}\h)$. 
    This Poisson $\SL(2,\K)$-algebra has dimension 18. 
    From these $\SL(2,\K)$-tables we deduce that the underlying Lie algebra of $H_E^{\bullet,\bullet}(\h)$ is isomorphic to  
$\Gl(3,\K)\ltimes\Gl(3,\K)_{ab}$ with the 
first factor acting on the second (abelian) one by the adjoint representation. 
We find it remarkable that the Lie algebra structure on 
$H_{E}^{\bullet,\bullet}(\h)$ contains a semisimple 
Lie subalgebra (in this case $\Sl(3,\K)$) 
strictly larger than the Levi factor 
of $\text{Der}(\h)$, which in this case is $\Sl(2,\K)\subset 
H^{1}(\h,\h)$.  This means that the Levi factor of the 
Lie algebra $H_{E}^{\bullet,\bullet}(\h)$ has nontrivial elements outside $H^{1}(\h,\h)$. 
Finally, this leads us to find a family of commutative 
Poisson algebras whose underlying Lie structure is $\Gl(n,\K)\ltimes\Gl(n,\K)_{ab}$ (arbitrary $n$) such that, for $n=3$, is isomorphic to 
$H_E^{\bullet,\bullet}(\h)$.
\end{abstract}

\maketitle

\section{Introduction}
This work is part of a project originally intended to investigate in detail the multiplicative structures that appear in various cohomology groups, particularly the Gerstenhaber and Poisson (super-)algebras associated with certain (co)homology spaces of associative or Lie algebras. 
It is widely acknowledged that obtaining a description of these types of algebras is a very challenging problem, even if we were satisfied
with some superficial information about the product in a small portion of the entire cohomology space.

In this context, the purpose of this article is twofold.
Firstly, we introduce the concept of the $G$-table of a given $G$-algebra (an algebra on which $G$ acts semisimply by algebra automorphisms, see \cite{AT95,JT95}).
If $\A$ is a $G$-algebra, then, roughly speaking, we have a $G$-module decomposition into irreducible 
$G$-submodules
\[
\A=\bigoplus_{r\in R} \A_r,
\]
and, in very vague terms, the
$G$-table of $\A$ is certain set of scalars
$c_{r_1,r_2}^{s}\in \K$ satisfying, in some sense,
\begin{equation}\label{eq.1}
\A_{r_1}\cdot \A_{r_2}=\sum_{s\in R} c_{r_1,r_2}^{s}\A_{s}
\end{equation}
resembling the ordinary multiplication table of $\A$
with respect to a given basis.
In \S\ref{sec.G_tables} we formalize this concept by providing a precise definition of a $G$-table of a given $G$-(co)algebra.
For instance, if $\A=M_3(\R)$ then $\A$ is a 
$SO(3)$-algebra and 
$\A=\A_0\oplus \A_1\oplus \A_2$ as a $SO(3)$-module, where $\A_1$ and $\A_2$ are, respectively, the skew-symmetric and symmetric traceless matrices, and $\A_0$ are the multiples of the identity.
The $SO(3)$-table of $\A$ turns out to be 
\begin{center}
\begin{tabular}{|c||c|c|c|}
\hline 
 & $\A_{0}$ & $\A_1$ & $\A_2$ \\
\hline 
\hline 
$\A_{0}$ & $\A_0$ &  $\A_1$ & $\A_2$ \rule[-2mm]{0mm}{7mm} \\
\hline 
$\A_{1}$ & $\A_1$ & $\A_0 + \A_1 + \A_{2}$ & $\A_1 + \A_2$ 
\rule[-2mm]{0mm}{7mm} \\
\hline 
$\A_{2}$ & $\A_2$ & $\A_1 - \A_{2}$ & $\A_0 + \A_1 + \A_{2}$ 
\rule[-2mm]{0mm}{7mm} \\
\hline 
\end{tabular}
\medskip

{The $SO(3)$-table of the associative algebra $M_3(\R)$.}
\end{center}
When computing a $G$-table, there are some choices that allow certain normalizations, this is why we were able to have, in this particular case, all the coefficients equal to $\pm 1$ in the previous $SO(3)$-table.  

The $G$-table of $\A$ provides very precise information about both, the algebra structure and the $G$-module structure of $\A$.
If $G=\{1\}$ is the trivial group, the $G$-table of an algebra is (almost) the same as the ordinary multiplication table of $\A$. 
On the other hand, when $G$ is large, the $G$-table is much more concise than 
an ordinary multiplication table and, still, 
the ordinary table can be deduced from the $G$-table.
When $G$ is `sufficiently' large, many fine details about the algebra structure of $\A$ are revealed by $G$-table. For instance, the above table shows us immediately that $M_3(\R)$ does not have any proper $SO(3)$-subalgebra (subalgebra and $SO(3)$-submodule). 
We can also observe that $[\A_1,\A_2]=2\A_2$ (which is the same as saying that $\A_2$ is a $\mathfrak{so}(3)$-submodule, note that $\A_1$ is the vector space $\mathfrak{so}(3)$). 
The first part of the paper includes many examples of $G$-tables and concludes by associating a plain algebra to a given $G$-algebra and analyzing some basic functoriality properties of this assignment.

\medskip

Obtaining the actual $G$-table of a given $G$-algebra $\A$ requires a considerable amount of work, in addition to a very precise set up in advance, 
that we refer to as the labeling of $\Rep_{G}$. 
But, it is worth it since it becomes a very powerful tool as it 
shown in the second part of the paper. 
In \S\ref{sec.Poisson_Struct} we describe the 
 Poisson algebra structure of
    the even-degree part of the cohomology associated to the cotangent bundle of the 3-dimensional Heisenberg Lie group. 
    Recall that, by a theorem of Nomizu \cite{N54}, the de Rham cohomology of 
the compact quotients of this cotangent bundle 
is isomorphic to the Lie algebra cohomology of $\h\ltimes\h^*$
where $\h$ is the 3-dimensional Heisenberg Lie algebra. 
We know that $\h$ is a Lie $\SL(2,\K)$-algebra and thus
the even cohomology 
\[
H_E^{\bullet,\bullet}(\h)=\bigoplus_{k \text{ even}}\;\bigoplus_{i+j=k} 
H^{i}(\h,\bwedge^j\h)
\]
is a Poisson $\SL(2,\K)$-algebra (of dimension 18).
First, we give the $\SL(2,\K)$-tables of the cup product and 
of the Poisson bracket in $H_E^{\bullet,\bullet}(\h)$. 
This requires three steps: (i) choosing a  labeling of $\Rep_{\SL(2,\K)}$, 
(ii) obtaining the  $\SL(2,\K)$-module decomposition of 
$H_E^{\bullet,\bullet}(\h)$ and matching it with the labeling, 
and (iii) a lengthy and very careful work to compute the actual coefficients of the $G$-tables (cup product and Lie bracket). 

From these $\SL(2,\K)$-tables we derive that the underlying Lie algebra of $H_E^{\bullet,\bullet}$ is isomorphic to  
$\Gl(3,\K)\ltimes\Gl(3,\K)_{ab}$ with the 
first factor acting on the second (abelian) one by the adjoint representation. 
It turns out that, for arbitrary $n$, 
the semidirect product (Lie algebra) $\Gl(n,\K)\ltimes\Gl(n,\K)_{ab}$, with 
associative commutative product 
\begin{align*}
(a_0I_n+A_0,\;& a_1I_n+A_1)\cdot (b_0I_n+B_0,\;b_1I_n+B_1) \\[3mm]
=
\Big(& a_0b_0 I_n   \;\;+\;\; a_0B_0 + b_0A_0\;,\;\\
&\big(a_0b_1+a_1b_0+\text{tr}(A_0B_1+A_1B_0)-\tfrac2{n}\text{tr}(A_0B_0)\big)I_n \;\;+\;\;
 a_0B_1 + b_0A_1+ A_0B_0+B_0A_0
 \Big)
\end{align*}
(here $A_0,A_1,B_0,B_1\in \Sl(n,\K)$)
is a Poisson algebra, and $H_E^{\bullet,\bullet}(\h)$
is isomorphic to the member of this family corresponding to $n=3$. 
The proof of this automorphism follows, in an almost straight forward way,
by comparing the corresponding $\SL(2,\K)$-tables.

It is remarkable for us that the Lie algebra structure on 
$H_{E}^{\bullet,\bullet}(\h)$ contains a semisimple 
Lie subalgebra (in this case $\Sl(3,\K)$) 
strictly larger than the Levi factor 
of $\text{Der}(\h)$, which in this case is $\Sl(2,\K)\subset 
H^{1}(\h,\h)$.  This means that the Levi factor of the 
Lie algebra $H_{E}^{\bullet,\bullet}(\h)$ has nontrivial elements outside $H^{1}(\h,\h)$.

\medskip

\noindent
\textbf{Acknowledgement.} 
We are deeply thankful to Marco Farinati 
for the very fruitful conversations on Poisson algebras that were crucial for our work.

\section{\texorpdfstring{$G$}{G}-tables}\label{sec.G_tables}

Let $\A$ be an algebra (not necessarily associative) over a field of characteristic zero $\K$.
If $G$ is a group and $G$ acts on $\A$ by algebra automorphisms in such a way that $\A$ is a direct sum of finite dimensional irreducible $G$-submodules, we  will say  that $\A$  is a \emph{$G$-algebra}. 
If $\A$ is a coalgebra, with comultiplication $\Delta$,
such that its dual $\A^*$ (with multiplication $\Delta^*$) is a $G$-algebra, we will say  that $\A$  is a \emph{$G$-coalgebra}. 

Recall that the 
\emph{multiplication table} of $\A$, with respect to a given 
basis  $B=\left\{ e_{r}: r\in R\right\}$, is the set of scalars $c_{r_1r,_2}^{s}\in \K$ satisfying
\[
e_{r_1}e_{r_2}=\sum_{s\in R} c_{r_1,r_2}^{s}e_{s} 
\]
(for each pair $r_1,r_2\in R$, the scalar $c_{r_1,r_2}^{s}$ is zero for all but a finite number of $s\in R$).

 For a $G$-algebra $\A$, we will introduce 
 the concept of the $G$-table of $\A$, which will be 
 a sort of multiplication table of $\A$ with respect to a given 
$G$-module decomposition of $\A$ in irreducible $G$-submodules.
Very roughly speaking, assume that 
\begin{equation}\label{eq.Decomp1}
\A=\bigoplus_{r\in R} \A_r
\end{equation}
is a $G$-module decomposition of $\A$ with $\A_r$ an irreducible $G$-submodule of $\A$. We want to define 
the multiplication $G$-table of $\A$, with respect to \eqref{eq.Decomp1}, as a set of scalars
$c_{r_1,r_2}^{s}\in \K$ satisfying (in some sense to be explained  below) 
\begin{equation}\label{eq.G-table_informal}
\A_{r_1}\cdot \A_{r_2}=\sum_{s\in R} c_{r_1,r_2}^{s}\A_{s}
\end{equation}
(as earlier, for each pair $r_1,r_2\in R$, the scalar $c_{r_1,r_2}^{s}$ will be zero for all but a finite number of $s\in R$).
In order to formalize this idea we need to introduce some preliminaries. 

\subsection{Labeling of \texorpdfstring{$\Rep(G)$}{Rep(G)}}\label{subsec.Lableling}
In what follows, $G$ will be a group
such that every finite dimensional representation of $G$ over 
$\K$ is completely reducible. We are mainly thinking in a semisimple 
algebraic group $G$ over $\K$. We will denote by 
$\Rep(G)$ the (tensor) category of finite dimensional representations of $G$.

In this context, we fix a set of representatives
\[
\hat{G}=\left\{ (\pi_i,V_{i}): i\in \hat I \right\}
\]
of  the isomorphism classes of irreducible representations in $\Rep(G)$.
In addition, given $i_1,i_2,j \in \hat I$, we also fix a basis 
$\bas^{i_1,i_2,j }$ of $\Hom_{G}\left(V_{i_1}\otimes V_{i_2},V_{j}\right)$, that is we choose $G$-morphisms, 
\[
 m^{i_1,i_2,j}_{q}:V_{i_1}\otimes V_{i_2}\to V_{j},\qquad q=1,\dots,d^{i_1,i_2,j},
 \]
where 
$d^{i_1,i_2,j}=\dim\big( \Hom_{G}\left(V_{i_1}\otimes V_{i_2},V_{j}\right)\big)$,
and we set 
\[
\bas^{i_1,i_2,j}=
\left\{ m^{i_1,i_2,j}_{q}: q=1,\dots,d^{i_1,i_2,j}
\right\}.
\]
\begin{Def}
    We say that the choice of the set $\hat{G}$ of representatives and the choice of each basis 
$\bas^{i_1,i_2,j}$ for $i_1,i_2,j\in \hat I$, 
constitute a \emph{labeling} of  $\Rep(G)$.
Occasionally, it might be enough, for our proposes, to consider only a subset
$\hat G_0$ of $\hat G$ (corresponding to $\hat I_0\subset\hat I$) and to choose  $\bas^{i_1,i_2,j }$ 
only for some $i_1,i_2,j \in \hat I_0$.
We will say that this choice constitute a \emph{partial labeling} of $\Rep(G)$.
\end{Def}

Note that giving a labeling of $\Rep(G)$
is almost the same as fixing a set of representatives
of $\hat{G}$ and fixing a particular choice of the 
 $3j$-symbols for the tensor category $\Rep(G)$.

\begin{remark}\label{rmk.abuse:language} 
In what follows we will frequently abuse
the notation identifying the index set $\hat I$ 
(resp. $\hat I_0$) with 
the set $\hat G$ (resp. $\hat G_0$).
\end{remark}

\begin{remark}\label{rmk.absense_q} 
For some groups $G$ it may happen that 
$\Hom_{G}\left(V_{i_1}\otimes V_{i_2},V_{j}\right)$ is 1-dimensional for all $i_1,i_2,j\in \hat I_0$;
in these cases we will omit the subindex $q$
running from $1$ to $d^{i_1,i_2,j}=1$.
A particular instance of this situation is when $G$ 
is a simply reducible group, see for instance \cite{kazarin2008finite} or \cite{hamermesh2012group}.
\end{remark}

\begin{remark}\label{rmk.enough_label-c} Although it is not mandatory, in many instances it 
is convenient to choose the bases 
$\bas^{i_1,i_2,j}$ so that, for all $q$, the $G$-morphisms
${m}^{i_1,i_2,j}_{q}$ and  
${m}^{i_2,i_1,j}_{q}$, $i_1\ne i_2$, coincide up to the canonical braiding
$V_{i_1}\otimes V_{i_2}\to V_{i_2}\otimes V_{i_1}$. 
Also, assuming that $i=0$ corresponds to the trivial representation of $G$, it is also convenient to choose the bases 
$\bas^{0,i,i}=\{{m}^{0,i,i}\}$ and  
$\bas^{i,0,i}=\{{m}^{i,0,i}\}$ with 
\[
{m}^{0,i,i}(x\otimes v)=
{m}^{i,0,i}(v\otimes x)=xv
\]
for all $x\in V_0=\K$ and $v\in V_i$.
Note that, here, we have omitted the index $q=1$ since 
$d^{0,i,i}=d^{i,0,i}=1$ for all $i\in \hat I$.
\end{remark}

Next, we give some examples that will be used in the following sections.

\begin{Ej}\label{ex.labeling_S_3}
\textbf{A labeling of $\Rep(\Sy_3)$.}
Let $G=\Sy_{3}=\{(),(12),(23),(13),(123),(132)\}$ be the permutation group on $3$ letters and 
let us label $\Rep(G)$. 
We have 
\[
\hat \Sy_{3}=\{(\pi_{tr},\K_{tr}),(\pi_{sg},\K_{sg}),(\pi_{std},\K^2_{std})\}
\]
 where $tr$ and $sg$ are, respectively, the trivial and sign representation,
and $std$ is the standard representation defined on $\K^2_{std}$ by
the following matrices given 
in the canonical basis $\{e_1,e_2\}$:
\[
\begin{matrix}
x=& () & (12) & (23) & (13) & (123) & (132) \\[2mm]
\pi_{std}(x)=&  \begin{pmatrix}1 & 0\\0 & 1\end{pmatrix} & 
                 \begin{pmatrix}0 & 1\\1 & 0\end{pmatrix} & 
                 \begin{pmatrix}1 & 0\\-1 & -1\end{pmatrix} & 
                 \begin{pmatrix}-1 & -1\\0 & 1\end{pmatrix} & 
                 \begin{pmatrix}-1 & -1\\1 & 0\end{pmatrix} & 
                 \begin{pmatrix}0 & 1\\-1 & -1\end{pmatrix}
\end{matrix}
\]
To complete the labeling of $\Rep(G)$ we need to give the (non-empty) bases 
$\bas^{i_1,i_2,j}$, with $i_1,i_2,j\in \hat \Sy_{3}$
(see Remark \ref{rmk.abuse:language}).
These are given as follows.

\begin{itemize}
\item We follow the convention explained  
in  Remark \ref{rmk.enough_label-c} which describes, in this case, the bases
\[
\bas^{tr,j,j}=\{m^{tr,j,j}\},\quad
\bas^{j,tr,j}=\{m^{j,tr,j}\},\quad\text{for $j=tr,sg,std$}.
\]

\item $\bas^{sg,sg,tr}=\{m^{sg,sg,tr}\}$ where 
$m^{sg,sg,tr}:{\K_{sg}}\otimes {\K_{sg}}\to \K_{tr}$ sends $1_{\K_{sg}}\otimes1_{\K_{sg}}\mapsto 1_{\K_{tr}}$.
\end{itemize}
In the following bases, the matrices are associated to the `lexicographic' basis of the corresponding tensor product. 
In particular, the ordered 
basis of $\K^2_{std}\otimes \K^2_{std}$ is 
$\{e_1\otimes e_1,e_1\otimes e_2,e_2\otimes e_1,e_2\otimes e_2\}$.
\begin{itemize}
\item $\bas^{sg,std,std}=\bas^{std,sg,std}=\left\{ \begin{pmatrix}1 & 2\\-2 & -1\end{pmatrix}\right\}$. 
This means  
\begin{equation}
\begin{split}
m^{sg,std,std}(1_{\K_{sg}}\otimes e_1)=m^{std,sg,std}(e_1\otimes 1_{\K_{sg}} )&=e_1-2e_2, \\[1mm]
m^{sg,std,std}(1_{\K_{sg}}\otimes e_2)=m^{std,sg,std}(e_2\otimes 1_{\K_{sg}} )&=2e_1-e_2.
\end{split}
\end{equation}
We emphasize that $m^{sg,std,std}$ and $m^{std,sg,std}$ are uniquely determined, up to scalar, since they are $G$-morphisms.

\item $\bas^{std,std,tr}=\left\{ \begin{pmatrix}2 & 1 & 1 & 2\end{pmatrix}\right\}$.
This means  
\begin{equation}
\begin{split}
m^{std,std,tr}(e_1\otimes e_1)=
m^{std,std,tr}(e_2\otimes e_2)&=2\,1_{\K_{tr}}, \\[1mm]
m^{std,std,tr}(e_1\otimes e_2)=
m^{std,std,tr}(e_2\otimes e_1)&=1_{\K_{tr}}.
\end{split}
\end{equation}

\item $\bas^{std,std,sg}=\left\{ \begin{pmatrix}0 & 1 & -1 & 0\end{pmatrix}\right\}$.
This means  
\begin{equation}
\begin{split}
m^{std,std,sg}(e_1\otimes e_1)&=
m^{std,std,sg}(e_2\otimes e_2)=0, \\[1mm]
m^{std,std,sg}(e_1\otimes e_2)&=
-m^{std,std,sg}(e_2\otimes e_1)=1_{\K_{sg}}.
\end{split}
\end{equation}

\item $\bas^{std,std,std}=\left\{ \begin{pmatrix}-1 & 1 & 1 & 2\\2 & 1 & 1 & -1\end{pmatrix}\right\}$.
This means  
\begin{equation}
\begin{split}
m^{std,std,std}(e_1\otimes e_1)&=-e_1+2e_2,\\[1mm]
m^{std,std,std}(e_1\otimes e_2)&=
m^{std,std,std}(e_2\otimes e_1)=e_1+e_2,\\[1mm]
m^{std,std,std}(e_2\otimes e_2)&=2e_1-e_2.
\end{split}
\end{equation}
\end{itemize}
This completes the labeling of $\Rep(\Sy_{3})$. 
Note that, in this example, $d^{i_1,i_2,j}\le 1$ for all $i_1,i_2,j\in \hat I$ and thus we omitted the subindex $q$
as indicated in Remark \ref{rmk.absense_q}.
\end{Ej}

\begin{Ej}\label{ex.G=GL_k}
\textbf{A partial labeling of $\Rep(GL(k,\K))$.}
 Let $G=GL(k,\K)$, in this example we will partially label  $\Rep(G)$ 
with 
\[
\hat G_0=\{(\pi_{0},\K_{0}),(\pi_{Ad},\Sl(k,\K))\}\subset \hat G,
\]
where $\pi_{0}$ and $\pi_{Ad}$ are, respectively, the trivial and the adjoint representation of $G$. 
As we did for $\Sy_3$, we choose 
\[
\bas^{0,0,0}=
\left\{ m^{0,0,0}\right\},
\qquad
\bas^{0,Ad,Ad}=
\left\{ m^{0,Ad,Ad}\right\},
\qquad
\bas^{Ad,0,Ad}=
\left\{ m^{Ad,0,Ad}\right\},
\]
as indicated by Remark \ref{rmk.enough_label-c}. 
In addition let 
\[
\bas^{Ad,Ad,0}=
\left\{ m^{Ad,Ad,0}\right\},
\qquad
\bas^{Ad,Ad,Ad}=
\left\{ m^{Ad,Ad,Ad}_1, m^{Ad,Ad,Ad}_2\right\}
\]
where
\begin{align*}
 m^{Ad,Ad,0}(A\otimes B)   & =\Traza(AB),  \\
 m^{Ad,Ad,Ad}_1(A\otimes B) & =[A,B]=AB-BA, \\
 m^{Ad,Ad,Ad}_2(A\otimes B) & =AB+BA-\frac{2}{k}\Traza(AB)\,I_k, 
\end{align*}
(here $I_k$ is the identity matrix). This completes the partial labeling of $\Rep(G)$.

Note that, if $k=2$, it is easy to see that $AB+BA-\frac{2}{k}\Traza(AB)\,I_k =0$ for all $A,B\in \Sl(2,\K)$. Therefore, $m^{Ad,Ad,Ad}_2 = 0$ when $k=2$ and $d^{Ad,Ad,Ad}=1$. On other hand, $d^{Ad,Ad,Ad}=2$ for all $k> 2$ and we have chosen the basis $\bas^{Ad,Ad,Ad}=\{m^{Ad,Ad,Ad}_1, m^{Ad,Ad,Ad}_2\}$ in such a way that their elements are, respectively, skew-symmetric and symmetric.
\end{Ej}

\medskip

In the following two examples we are going to consider the group $G=SL(2,\K)$. A standard labeling of $\Rep(G)$ is given by the classical Wigner $3j$-symbols (or Clebsh-Gordan coefficients, see for instance \cite{carter1995classical}).
We will not need in this article a complete labeling of $\Rep(G)$
and therefore, instead of recalling the classical  Clebsh-Gordan
coefficients we will present two partial labelings of $\Rep(G)$
that will suffice for our proposes. 

\begin{Ej}\label{ex.6J-labeling} 
\textbf{First partial labeling of $\Rep(SL(2,\K))$.} 
Let $G=SL(2,\K)$, in this example we enlarge a little bit, for $k=2$, 
the partial labeling of $\Rep(GL(k,\K))$
given in Example \ref{ex.G=GL_k}.
Let 
\[
\hat G_0=\{(\pi_{0},\K_{0}),(\pi_{1},\K^2),(\pi_{2},\Sl(2,\K))\}\subset \hat G,
\]
here, the subindex in $\pi$ represents the highest weight of the representation. 
This means that 
$\pi_{0}$ is the trivial representation, 
$\pi_{1}$ is the canonical representation, and 
$\pi_{2}$ is the adjoint representation. 
We follow Remark \ref{rmk.enough_label-c}
to choose 
\[
\bas^{0,i,i}=
\left\{ m^{0,i,i}\right\},
\quad
\bas^{i,0,i}=
\left\{ m^{i,0,i}\right\},\qquad\text{for $i=0,1,2$.}
\]
In addition let 
\[
\bas^{1,1,0}=
\left\{m^{1,1,0}\right\},
\quad
\bas^{1,1,2}=
\left\{m^{1,1,2}\right\},
\quad
\bas^{2,1,1}=
\left\{m^{2,1,1}\right\},
\quad
\bas^{2,2,0}=
\left\{m^{2,2,0}\right\},
\quad
\bas^{2,2,2}=
\left\{m^{2,2,2}\right\}
\]
where
\begin{align*}
 m^{1,1,0}\big((x_1,x_2)\otimes (y_1,y_2)\big) &= \det\begin{pmatrix}
 x_1&x_2\\y_1&y_2\end{pmatrix} \\
 m^{1,1,2}\big((x_1,x_2)\otimes (y_1,y_2)\big) &=
 \begin{pmatrix}
 x_1y_2+x_2y_1 &-2x_1y_1\\2x_2y_2&-x_1y_2-x_2y_1 \end{pmatrix} \\
 m^{2,1,1}\big(A\otimes(x_1,x_2)\big) &= (x_1,x_2).\,A^t \\
 m^{2,2,0}(A\otimes B) &= \Traza(AB),  \\
 m^{2,2,2}(A\otimes B) &= [A,B]=AB-BA.
\end{align*}
Following Remark \ref{rmk.enough_label-c}, we also choose 
$m^{1,2,1}\big((x_1,x_2)\otimes A\big)=m^{2,1,1}\big(A\otimes(x_1,x_2)\big)$.

For later use, we pick for each representation in $\hat G_0$
the following highest weight vector $\bar v$ with respect to the standard
$\mathfrak{s}$-triple of $\mathfrak{sl}(2,\K)$: 
\begin{itemize}
\item For $(\pi_{0},\K_{0})$, let $\bar v=1_{\K_{0}}$.
\item For $(\pi_{1},\K^2)$, let $\bar v=(1,0)$.
\item For $(\pi_{2},\Sl(2,\K))$, let 
$\bar v=\left(\begin{smallmatrix}0&1\\0&0\end{smallmatrix}\right)$.
\end{itemize}
\end{Ej}

\begin{Ej}\label{ex.partial_G=SL_2}
\textbf{Second partial labeling of $\Rep(SL(2,\K))$.} 
Again, let $G=SL(2,\K)$, but now, $\hat G_0=\hat G$. 
This will be a partial labeling of $\Rep(SL(2,\K))$ becuse we will not make a choice of 
$\bas^{i_1,i_2,j}$ for all $i_1,i_2,j\in I$.

It is well known that 
$\K[x,y]$ contains all 
finite dimensional irreducible $G$-modules exactly once.
Here, $g\in SL(2,\K)$ acts on $p\in \K[x,y]$ by $(g.p)(x,y)=p((x,y)g)$.
If 
$\K[x,y]_r$, $r\in \Z_{\ge 0}$,
is the subspace of $\K[x,y]$ consisting of all homogeneous polynomials of degree $r$, then $\K[x,y]_r$ is irreducible as $SL(2,\K)$-module.
Thus we take $\hat G=\{\K[x,y]_r: r\in \Z_{\ge 0}\}$ and we choose the bases $\bas^{r_1,r_2,r_1+r_2}=\{m^{r_1,r_2,r_1+r_2} \}$, by using the polynomial multiplication,
that is
\[
m^{r_1,r_2,r_1+r_2}(p_1,p_2)=p_1p_2,\qquad p_i\in \K[x,y]_{r_i}.
\]
This is a partial labeling of $\Rep(G)$ since, in order to complete it, we would need to define
$m^{r_1,r_2,s}$ for all $s\in\Z$ such that  $|r_1-r_2|\le s \le r_1+r_2$ and  $s\equiv r_1+r_2\mod (2)$. 
\end{Ej}

\subsection{\texorpdfstring{$G$}{G}-tables of algebras (and coalgebras)}\label{subsec.G-tables}
Once $\Rep(G)$ has been (partially) labeled, we are 
almost ready to define (or compute)  the 
$G$-table of 
any $G$-algebra $\A$. 
First, from the definition of a $G$-algebra, we can fix a decomposition 
\[
\A=\bigoplus_{r\in R}\A_{r}.
\]
into irreducible $G$-submodules.
Here $\A_{r}$ is finite dimensional and isomorphic to some $V_i\in \hat{G}$ for some $i=\ep(r)\in \hat I$. 
Thus, we may choose a $G$-morphism 
\[
\tau_{r}:V_{\ep(r)}\rightarrow\A
\] 
whose image is $\A_{r}$. We denote by $\bar\tau_{r}:V_{\ep(r)}\to\A_{r}$ the corresponding $G$-isomorphism. 

Now, for each pair $r_1,r_2\in R$,  we consider the following basis of $\Hom_{G}\left(\A_{r_1}\otimes\A_{r_2},\A_{s}\right)$ obtained from the labeling of $\Rep(G)$
\[
\tilde{\bas}^{r_1,r_2,s}=\left\{ 
\tilde{m}^{r_1,r_2,s}_{q}:q=1,\dots, d^{\ep(r_1),\ep(r_2),\ep(s)}\right\},
\]
where, for each $q=1,\dots,d^{\ep(r_1),\ep(r_2),\ep(s)}$, 
\[
\tilde{m}^{r_1,r_2,s}_{q}=
\bar\tau_{s}\circ
m^{\ep\left(r_1\right),\ep\left(r_2\right),\ep\left(s\right)}_{q}\circ
\left(\bar\tau_{r_1}\otimes\bar\tau_{r_2}\right)^{-1}.
\]
Let $\mu:\A\otimes\A\to\A$ denote the multiplication of $\A$. 
Since, for each pair $r_1,r_2\in R$, the restriction 
\[
\mu_{r_1,r_2}:\A_{r_1}\otimes\A_{r_2}\to\A
\]
is a $G$-morphism, 
we can express it as a linear combination of the elements in
$\bigcup_{s\in R}\tilde{\bas}^{r_1,r_2,s}$, that is 
\begin{equation}\label{eq: G-table}
   \mu_{r_1,r_2}=\sum_{s\in R}\sum_{q=1}^{d^{\ep(r_1),\ep(r_2),\ep(s)}}c_{r_1,r_2}^{s,q}\;\tilde{m}^{r_1,r_2,s}_{q}
\end{equation}
for some uniquely determined scalars 
$c_{r_1,r_2}^{s,q}\in\K$. 
Since $\A_{r_1}\otimes\A_{r_2}$ is finite dimensional, 
it follows that, for each pair $r_1,r_2\in R$,
only for a finite number of $s\in R$ there are some 
$q\ge 1$ such that $c_{r_1,r_2}^{s,q}\ne0$.

Summarizing, we give the following definition.

\begin{Def} \label{Def: G-table} We say that the set
of scalars  $c_{r_1,r_2}^{s,q}\in\K$, defined by equation \eqref{eq: G-table}, is the \emph{$G$-table} of $\A$. 
We point out that a $G$-table depends on
\begin{enumerate}
\item The labeling of $\Rep(G)$, that is the choice of $\hat{G}$ and the choice of each basis 
$\bas^{i_1,i_2,j}$ for $i_1,i_2,j\in \hat I$. 

\item The choice of the $G$-morphisms 
$\tau_{r}:V_{\ep(r)}\rightarrow\A$, $r\in R$, so that 
$\A=\bigoplus_{r\in R} \tau_{r}(V_{\ep(r)})$.
\end{enumerate}
If $\A$ is a $G$-coalgebra, its \emph{$G$-cotable} is the 
$G$-table of $\A^*$.
\end{Def}

The $G$-table 
\eqref{eq: G-table} of $\A$ 
is a very synthetic way to explicitly 
describe the multiplicative structure of $\A$ in combination with 
its $G$-module structure. 
It contains much more information than an ordinary
multiplication table. Indeed, since we labeled 
$\Rep(G)$, we have a an explicit description of
each $m^{i,j,k}_q$, and thus we  
have the multiplication table of 
each $m^{i,j,k}_q$, therefore it is straightforward to give the 
ordinary multiplication table of $\A$ from its $G$-table 
\eqref{eq: G-table}.

\medskip

We collect now some observations that will be frequently used in the rest or the paper. 

\begin{remark}\label{rmk.enough_label-a}
If we want to compute a 
$G$-table of a single $G$-algebra $\A$, it is enough to 
partially label $\Rep(G)$ with $\hat G_0$ equal to the 
subset of $\hat G$ consisting of the irreducible $G$-modules 
appearing in the  $G$-module decomposition of $\A$.
\end{remark}

\begin{remark}\label{rmk.notation_table}
When some $d^{i_1,i_2,j}$ are greater than 1, 
\eqref{eq: G-table} looks a bit more complicated than the 
announced \eqref{eq.G-table_informal}. 
Nevertheless, in order to keep in mind the idea that the $G$-table 
is a generalization of the ordinary multiplication table
(see Example \ref{ex.G=1} below), 
we will frequently write 
\begin{equation}\label{eq: G-table_with_A}
    \A_{r_1}\cdot \A_{r_2}=\sum_{s\in R}\sum_{q=1}^{d^{\ep(r_1),\ep(r_2),\ep(s)}}c_{r_1,r_2}^{s,q}\;\A_{s,q}
\end{equation}
instead of \eqref{eq: G-table}. We insist that this is only a notation that describes \eqref{eq: G-table}.
We will also depict the $G$-table \eqref{eq: G-table_with_A} in a table with rows and columns corresponding to $\A_{r_1}$ and $\A_{r_2}$ respectively, and the associated entry will be the right hand side of  \eqref{eq: G-table_with_A}
(see, for instance, Example \ref{ex.G-table_A=Mkk} below).
\end{remark}

\medskip

We now give some basic examples of $G$-tables.

\begin{Ej}
\label{ex.G=1} Ordinary multiplication tables of an algebra $\A$ are the 
$G$-tables corresponding to the trivial group $G=1$.
Indeed, $\hat G$  consists only 
of the trivial module $V_1=\K$, 
$\Hom_{G}\left(V_{1}\otimes V_{1},V_{1}\right)$ is 1-dimensional and we can choose 
$m^{1,1,1}$ as the map defined by $1_{\K}\otimes 1_{\K}\mapsto 1_{\K}$.
Choosing a basis 
$B=\left\{ e_{r}: r\in R\right\}$ of $\A$
is almost the same as choosing a decomposition of $\A$
into irreducible $G$-submodules, once $B$ has been chosen, the morphisms $\bar\tau_r:\K\to\K e_r\subset\A$ are determined up to a scalar and if we choose them so that $1_{\K}$ is 
mapped to $e_r$, then the associated $G$-table
is same as the multiplication table of $\A$ with respect to the basis $B$.
\end{Ej}
 
\begin{Ej} \label{ex.G-table_A=Mkk} Let $\A=M_{k}(\K)$ be the (associative) algebra of  $k\times k$ matrices over $\K$ and let $G=GL(k,\K)$ act on $\A$ 
by conjugation. The decomposition of $\A$ into irreducible
$G$-submodules is 
\[
\A = \A_0 \oplus \A_{1}
\]
where $\A_0=\K I_k$ ($I_k$ is the identity matrix) 
and $\A_1$ is subspace of matrices $A$ with $\Traza(A)=0$. 
In Example \ref{ex.G=GL_k} we partially labeled $\Rep(G)$ as suggested in Remark  \ref{rmk.enough_label-a}.
In particular, since is $I_k$ is the identity of $\A$ it follows that 
\[
\mu_{0,0}=
\tilde{m}^{0,0,0},\qquad
\mu_{0,1}=
\tilde{m}^{0,1,1},\qquad
\mu_{1,0}=
\tilde{m}^{1,0,1}.
\]
Since 
$AB=\frac{1}{k}\Traza(AB)I_k+\frac12 [A,B]+\frac12 (AB+BA-\frac{2}{k}\Traza(AB)I_k)$
it follows that 
\[
 \mu_{1,1}=\frac{1}{k}\tilde{m}^{1,1,0}+\frac12 \tilde{m}^{1,1,1}_1+\frac12 \tilde{m}^{1,1,1}_2.
\]
Thus, the $G$-table
of $M_{k}(\K)$ (with respect to the partial labeling given in Example \ref{ex.G=GL_k})  is
\begin{center}
\begin{tabular}{|c||c|c|}
\hline 
 & $\A_{0}$ & $\A_1$ \\
\hline 
\hline 
$\A_{0}$ & $\tilde{m}^{0,0,0}$ &  $\tilde{m}^{0,1,0}$ \rule[-2mm]{0mm}{7mm} \\
\hline 
$\A_{1}$ & $\tilde{m}^{1,0,0}$ & $\frac{1}{k} \tilde{m}^{1,1,0} +
\frac12 \tilde{m}^{1,1,1}_1 +\frac12 \tilde{m}^{1,1,1}_2$ 
\rule[-2mm]{0mm}{7mm} \\
\hline 
\end{tabular}
\end{center}
Following Remark \ref{rmk.notation_table}, we prefer to 
describe the $G$-table of $M_{k}(\K)$ as
\begin{center}
\begin{tabular}{|c||c|c|}
\hline 
 & $\A_{0}$ & $\A_1$ \\
\hline 
\hline 
$\A_{0}$ & $\A_0$ &  $\A_1$ \rule[-2mm]{0mm}{7mm} \\
\hline 
$\A_{1}$ & $\A_1$ & $\frac{1}{k} \A_0 +
\frac12 \A_{1,1} +\frac12 \A_{1,2}$ 
\rule[-2mm]{0mm}{7mm} \\
\hline 
\end{tabular}
\end{center}
(note the additional subindeces $q$ in $\A_{1,q}$).
We warn that, at first glance, this $G$-table may let the reader to mistakenly think that $\A$ is a commutative algebra. 
The non-commutativity of $\A$ is revealed by the non trivial linear combination of the skew-symmetric  $\tilde{m}^{1,1,1}_1$
and the symmetric maps $\tilde{m}^{1,1,0}$ and $\tilde{m}^{1,1,1}_2$.

We can also use the partial labeling for $G=GL(k,\K)$
given in Example \ref{ex.G=GL_k} to give the $G$-table of the Lie algebra $\A=\Gl(k,\K)$. It is straightforward to see that it is 
\begin{center}
\begin{tabular}{|c||c|c|}
\hline 
 & $\A_{0}$ & $\A_1$ \\
\hline 
\hline 
$\A_0$ & $0$ &  $0$ \rule[-2mm]{0mm}{7mm} \\
\hline 
$\A_1$ & $0$ & $ \A_{1,1}$ 
\rule[-2mm]{0mm}{7mm} \\
\hline 
\end{tabular}
\end{center}
\end{Ej}

 \medskip
 
\begin{Ej}  \label{ex.A=K[x,y]}
Let $G=\SL(2,\K)$. 
We can use the partial labeling of $\Rep(G)$ given in 
Example \ref{ex.partial_G=SL_2} to compute
the $G$-table of the $G$-algebra 
$\A=\K[x,y]$ (recall that $g\in SL(2,\K)$ acts on $p\in \A$ as
$(g.p)(x,y)=p((x,y)g)$).
Note that, even though $\A$ is not a finite dimensional $G$-module, it is isomorphic to the direct sum of finite dimensional irreducible  representations of $G$, with all of them appearing once.
It is straightforward to see that, with the labeling given 
in Example \ref{ex.partial_G=SL_2}, the $G$-table of $\A$ is 
\[
\A_{r_1}\cdot \A_{r_2}=
\A_{r_1+r_2}.
\]
\end{Ej}

 \medskip
 
\begin{Ej}  \label{ex.A=K[S3]}
Let $G=\Sy_3$ and let $\A=\K[\Sy_3]$ with $G$ acting 
on $\A$ by conjugation. In this case $\A$ is a Hopf algebra
and in this example we compute the $G$-table and the $G$-cotable
of $\A$ with respect to the labeling given in Example \ref{ex.labeling_S_3}.

The $G$-module decomposition of $\A$ is 
\[
\A=\A_{tr,1}\oplus \A_{tr,2}\oplus \A_{tr,3}\oplus 
\A_{sg} \oplus \A_{std}
\]
where, as $G$-modules (see Example \ref{ex.labeling_S_3})
\begin{align*}
\A_{tr,i}=\text{span}\{1_{i}\}&\simeq \K_{tr},\qquad  \tau_{tr,i}(1_{\K_{tr}})=1_i,\;i=1,2,3, \\
\A_{sg}=\text{span}\{s_{sg}\}&\simeq \K_{sg}, \qquad  \tau_{sg}(1_{\K_{sg}})=s_{sg}, \\
\A_{std} =\text{span}\{u_1,u_2\}&\simeq \K^2_{std}, 
\qquad  \tau_{std}(e_j)=u_j,\;j=1,2,
\end{align*}
and 
\begin{align*}
1_{1} & =\frac{1}{6}\big(()+(12)+(23)+(13)+(123)+(132)\big), 
& s_{sg} & =(123)-(132),\\
1_{2} & =\frac{1}{6}\big(()-(12)-(23)-(13)+(123)+(132)\big), 
& u_1 & =(12)-(23),\\
1_{3} & =\frac{1}{3}\big(2()-(123)-(132)\big), 
& u_2 & =(12)-(13).
\end{align*}
It follows that the $G$-table of $\A$ is (note that we replaced the 1-dimensional submodules by the given generator):
\begin{center}
\begin{tabular}{|c||c|c|c|c|c|}
\hline 
 & {$1_{1}$} & {$1_{2}$} & {$1_{3}$} & {$s_{sg}$} & {$\A_{std}$}\rule[-2mm]{0mm}{7mm} \\
\hline 
\hline 
{$1_{1}$} & {$1_{1}$} & {$0$} & {$0$} & {$0$} & {$0$}\rule[-2mm]{0mm}{7mm}  \\
\hline 
{$1_{2}$} & {$0$} & {$1_{2}$} & {$0$} & {$0$} & {$0$}\rule[-2mm]{0mm}{7mm}  \\
\hline 
{$1_{3}$} & {$0$} & {$0$} & {$1_{3}$} & {$s_{sg}$} & {$\A_{std}$}\rule[-2mm]{0mm}{7mm}  \\
\hline 
{$s_{sg}$} & {$0$} & {$0$} & {$s_{sg}$} & {$-3\,1_{3}$} & {$\A_{std}$}\rule[-2mm]{0mm}{7mm}  \\
\hline 
{$\A_{std}$} & {$0$} & {$0$} & {$\A_{std}$} & {$-\A_{std}$} & {$\frac{3}{2}\,1_{3}+\frac{3}{2}\,s_{sg}$}\rule[-2mm]{0mm}{7mm}  \\
\hline 
\end{tabular}
\end{center}
Some computations are needed to show that this is, indeed, the $G$-table.
For instance:
\begin{itemize}
\item the entry $s_{sg}\cdot s_{sg}= -3\,1_{3}$ follows from
\[
s_{sg}^2= -3\,1_{3} \quad\text{and}\quad
m^{sg,sg,tr}(1_{\K_{sg}},1_{\K_{sg}})=1_{\K_{tr}},
\] 
\item the entry $s_{sg}\cdot\A_{std}=\A_{std}$ follows from
\begin{align*}
s_{sg}u_1=-(12)-(23)+2(13)=u_1-2u_2\quad\text{and}\quad
m^{sg,std,std}(1_{\K_{sg}}\otimes e_1)=e_1-2e_2, \\
s_{sg}u_2=(12)-2(23)+(13)=2u_1-u_2\quad\text{and}\quad
m^{sg,std,std}(1_{\K_{sg}}\otimes e_2)=2e_1-e_2; 
\end{align*}
\item the entry $\A_{std}\cdot\A_{std}= \frac{3}{2}\,1_{3}+\frac{3}{2}\,s_{sg}$ follows from
\begin{align*}
u_1u_1=3\,1_{3}
&\quad\text{and}\quad
\begin{array}{l}
m^{std,std,tr}(e_1\otimes e_1)=2\,1_{\K_{tr}}, \\
m^{std,std,sg}(e_1\otimes e_1)=0, 
\end{array} \\[2mm]
u_1u_2=\tfrac32(1_{3}+s_{sg})
&\quad\text{and}\quad
\begin{array}{l}
m^{std,std,tr}(e_1\otimes e_2)=1_{\K_{tr}}, \\
m^{std,std,sg}(e_1\otimes e_2)=1_{\K_{sg}}, 
\end{array} \\[2mm]
u_2u_1=\tfrac32(1_{3}-s_{sg})
&\quad\text{and}\quad
\begin{array}{l}
m^{std,std,tr}(e_2\otimes e_1)=1_{\K_{tr}}, \\
m^{std,std,sg}(e_2\otimes e_1)=-1_{\K_{sg}}, 
\end{array} \\[2mm]
u_2u_2=3\,1_{3}
&\quad\text{and}\quad
\begin{array}{l}
m^{std,std,tr}(e_2\otimes e_2)=2\,1_{\K_{tr}}, \\
m^{std,std,sg}(e_2\otimes e_2)=0.
\end{array} 
\end{align*}
\end{itemize} 
Recall that the $G$-cotable of $\A$ is the $G$-table of $\A^*$ with product given by the transpose of
$\Delta:\A\rightarrow\A\otimes\A$,  
$\Delta(\sigma)=\sigma\otimes \sigma$, $\sigma\in G$.
Let us identify $\A^*$ with $\A$ via the inner product 
in which the basis $G$ of $\A$ is orthonormal. 
It follows that the product given by $\Delta$ is 
\[
g_1g_2=\begin{cases}
g_1,&\text{if $g_1=g_2$;} \\
0,&\text{if $g_1\ne g_2$.}
\end{cases}
\]
Now, the $G$-cotable of $\A$ is 
\begin{center}
\begin{tabular}{|c||c|c|c|c|c|}
\hline 
 & {$1_{1}$} & {$1_{2}$} & {$1_{3}$} & {$s_{sg}$} & {$\A_{std}$}\rule[-2mm]{0mm}{7mm} \\
\hline 
\hline 
{$1_{1}$}  & {$\tfrac{1}{6}\,1_{1}$} & {$\tfrac{1}{6}\,1_{2}$} & {$\tfrac{1}{6}\,1_{3}$} & {$\tfrac{1}{6}\,s_{sg}$} & {$\tfrac{1}{6}\,\A_{std}$}\rule[-2mm]{0mm}{7mm}  \\
\hline 
{$1_{2}$}  & {$\tfrac{1} {6}\,1_{2}$} & {$\tfrac{1} {6}\,1_{1}$} & {$\tfrac{1} {6}\,1_{3}$} & {$\tfrac{1} {6}\,s_{sg}$} & {$-\tfrac{1} {6}\,\A_{std}$}\rule[-2mm]{0mm}{7mm}  \\
\hline 
{$1_{3}$} & {$\tfrac{1} {6}\,1_{3}$} & {$\tfrac{1} {6}\,1_{3}$} & {$\frac23 1_1+\frac23 1_2+\frac13 1_3$} & {$-\frac13\, s_{sg}$} & {$0$}\rule[-2mm]{0mm}{7mm}  \\
\hline 
{$s_{sg}$} & {$\tfrac{1} {6}\,s_{sg}$} & {$\tfrac{1} {6}\,s_{sg}$} & {$-\frac13\, s_{sg}$} & {$2\, 1_1+2\, 1_2- 1_3$} & {$0$}\rule[-2mm]{0mm}{7mm}  \\
\hline 
{$\A_{std}$} & {$\tfrac{1} {6}\,\A_{std}$} & {$-\tfrac{1} {6}\,\A_{std}$} & {$0$} & {$0$} & {$1_1- 1_2+\frac13\, \A_{std}$}\rule[-2mm]{0mm}{7mm}  \\
\hline 
\end{tabular}
\end{center}
\end{Ej}

\medskip

\begin{Ej}  \label{ex.A=sl(3)}
Let $G=SL(2,\K)$ and $\A=\mathfrak{sl}(3,\K)$. 
We think $G$ as a subgroup  $SL(3,\K)$ located in the upper left corner. This way, $G$ acts on $\A$ by conjugation (it is the restriction to $SL(2,\K)$ of Adjoint action of $SL(3,\K)$ on $\mathfrak{sl}(3,\K)$). 
Thus, the $SL(2,\K)$-module decomposition of $\A$ is: 
\begin{align*}
\A&=
\K
\begin{pmatrix}
1 & 0 & 0 \\
0 & 1 & 0 \\
0 & 0 & -2 
\end{pmatrix}
\oplus
\begin{pmatrix}
* & * & 0 \\
* & -* & 0 \\
0 & 0 & 0 
\end{pmatrix}
\oplus
\begin{pmatrix}
0 & 0 & * \\
0 & 0 & * \\
0 & 0 & 0 
\end{pmatrix}
\oplus
\begin{pmatrix}
0 & 0 & 0 \\
0 & 0 & 0 \\
* & * & 0 
\end{pmatrix} \\
&\simeq
V_0\oplus V_2\oplus V_1\oplus V_1'
\end{align*} 
Considering the labeling of $\Rep(SL(2,\K))$ given in Example \ref{ex.6J-labeling}, the $SL(2,\K)$-table of $\mathfrak{sl}(3,\K)$ is
\begin{center}
\begin{tabular}{|c||c|c|c|c|}
\hline 
 & {$V_0$} & {$V_2$} & {$V_1$} & {$V_1'$} \rule[-2mm]{0mm}{6mm} \\
\hline 
\hline 
 {$V_0$} & & & {$3\,V_1$} & {$-3\,V_1'$} \rule[-1mm]{0mm}{6mm} \\
\hline 
 {$V_2$} &  & {$V_2$} & {$V_1$} & {$V_1'$} \rule[-1mm]{0mm}{6mm} \\
\hline 
 {$V_1$} & {$-3\,V_1$} & {$-V_1$} &  & {$-\frac{1}{2}\,V_2+\frac{1}{2}\,V_0$} \rule[-2mm]{0mm}{6mm} \\
\hline 
 {$V_1'$} & {$3\,V_1$}  & {$-V_1'$} & {$\frac{1}{2}\,V_2+\frac{1}{2}\,V_0$} &  \rule[-2mm]{0mm}{6mm} \\
\hline 
\end{tabular}
\end{center}
\end{Ej}

\subsection{\texorpdfstring{$G$}{G}-algebra morphisms, $G$-matrix, functoriality}
Let $G$ be a group satisfying the conditions required in \S\ref{subsec.Lableling} and assume that $\Rep(G)$ is labeled.
Given two $G$-algebras 
$\A$ and $\B$, a 
\emph{$G$-algebra morphism} between them is an algebra homomorphism
$\varphi:\A\to\B$ that, in addition, is a $G$-module morphism.

Let $\mu^{\A}$ and $\mu^{\B}$ denote, respectively, the multiplications 
of $\A$ and $\B$, and let 
$\A=\oplus_{r\in R}\A_r$ and 
$\B=\oplus_{x\in X}\B_x$ be $G$-module decompositions
with corresponding choices of 
$G$-module morphisms 
\begin{align*}
\tau^{\A}_{r}&:V_{\ep^{\A}(r)}\rightarrow\A_r\subset\A,\qquad
\epsilon^{\A}:R\to \hat I, \\
\tau^{\B}_{x}&:V_{\ep^{\B}(x)}\rightarrow\B_x\subset\B,\qquad
\epsilon^{\B}:X\to \hat I.
\end{align*}
For each $r\in R$ and each $x\in X$, 
let 
\begin{align*}
X_r&=\{y\in X:\epsilon^{\B}(y)= \epsilon^{\A}(r) \},\\ 
R_x&=\{s\in R:\epsilon^{\A}(s)= \epsilon^{\B}(x)\}.
\end{align*} 
Given a 
$G$-module morphism 
$\varphi:\A\to\B$, 
there exist scalars $f_{x,r}\in\K$, $r\in R$, $x\in X$,
 such that 
\begin{equation}\label{eq.f_rx}
\varphi_r=
\varphi|_{\A_r}=
\sum_{x\in X_r}
f_{x,r}\, \bar\tau^{\B}_{x}\circ(\bar\tau^{\A}_{r})^{-1}.
\end{equation}
Note that, since $\A_r$ is finite dimensional,
it follows that: (*) for  each $r\in R$, $f_{x,r}\ne 0$
only for a finite number of $x\in X_r$.
Conversely, given scalars $f_{x,r}\in\K$ satisfying
(*), it is clear that 
\eqref{eq.f_rx} defines a 
$G$-module morphism 
$\varphi:\A\to\B$. 
We call this set of scalar the \emph{$G$-matrix} of $f$ (with respect to the given decompositions of $\A$ and $\B$).
We want to describe, in terms of the $G$-tables of $\A$ and $\B$, 
the condition satisfied by the scalars $f_{x,r}\in\K$ when 
$\varphi$ is a $G$-algebra morphism.

Let
\begin{align*}
   \mu^{\A}_{r_1,r_2}&=\sum_{s\in R}\sum_{q=1}^{d^{\ep(r_1),\ep(r_2),\ep(s)}}c_{r_1,r_2}^{s,q}\;
\bar\tau^{\A}_{s}\circ
m^{\ep\left(r_1\right),\ep\left(r_2\right),\ep\left(s\right)}_{q}\circ
\left(\bar\tau^{\A}_{r_1}\otimes\bar\tau^{\A}_{r_2}\right)^{-1}, 
\quad c_{r_1,r_2}^{s,q}\in\K, \\
   \mu^{\B}_{x_1,x_2}&=\sum_{y\in X}\sum_{q=1}^{d^{\ep(x_1),\ep(x_2),\ep(y)}}d_{x_1,x_2}^{y,q}\;
\bar\tau^{\B}_{y}\circ
m^{\ep\left(x_1\right),\ep\left(x_2\right),\ep\left(y\right)}_{q}\circ
\left(\bar\tau^{\B}_{x_1}\otimes\bar\tau^{\B}_{x_2}\right)^{-1},
\quad d_{x_1,x_2}^{y,q}\in\K, 
\end{align*}
be the corresponding $G$-tables.
Note that, in order to simplify the notation, we have dropped the superscripts
$\A$ and $\B$ in both functions $\epsilon$. 
The difference between these two $\epsilon$ should be clear from the name of the variables ($r,s\in R$, $x,y\in X$).

\begin{thm}\label{thm.G-morphism}
Let  
$\varphi:\A\to\B$ be $G$-module morphism and let 
$f_{x,r}\in\K$ ($r\in R$, $x\in X_r$) be the $G$-matrix as defined in 
\eqref{eq.f_rx}.
Then, $\varphi$ is a $G$-algebra morphism if and only if 
\begin{equation}\label{eq.fxr_condition}
\sum_{s\in R_y}
c_{r_1,r_2}^{s,q}
\,
f_{y,s}
=
\sum_{x_1\in X_{r_1}}
\sum_{x_2\in X_{r_2}}
d_{x_1,x_2}^{y,q}\, 
f_{x_1,r_1}\, f_{x_2,r_2}
\end{equation}
for all $r_1,r_2\in R$, $y\in X$ and 
all $q=1,\dots,d^{\ep(r_1),\ep(r_2),\ep(y)}$.
\end{thm}
\begin{proof} The
$G$-module morphism 
$\varphi$ is a $G$-algebra morphism if and only if 
$\varphi\circ\mu^{\A}=
\mu^{\B}\circ (\varphi\otimes \varphi)$ and this means
\[
\varphi\circ\mu^{\A}_{r_1,r_2} =
\mu^{\B}\circ (\varphi_{r_1}\otimes \varphi_{r_2}),
\quad\text{for all }r_1,r_2\in R.
\]
On the one hand we have
\begin{align*}
\varphi\circ\mu^{\A}_{r_1,r_2} 
&=\sum_{s\in R}\sum_{q=1}^{d^{\ep(r_1),\ep(r_2),\ep(s)}}c_{r_1,r_2}^{s,q}\;
\varphi\circ
\bar\tau^{\A}_{s}\circ
m^{\ep\left(r_1\right),\ep\left(r_2\right),\ep\left(s\right)}_{q}\circ
\left(\bar\tau^{\A}_{r_1}\otimes\bar\tau^{\A}_{r_2}\right)^{-1} \\ 
&=\sum_{s\in R}
\sum_{q=1}^{d^{\ep(r_1),\ep(r_2),\ep(s)}}
\sum_{y\in X_s}
c_{r_1,r_2}^{s,q}
\;
f_{y,s}
\;
\bar\tau^{\B}_{y}\circ
m^{\ep\left(r_1\right),\ep\left(r_2\right),\ep\left(s\right)}_{q}\circ
\left(\bar\tau^{\A}_{r_1}\otimes\bar\tau^{\A}_{r_2}\right)^{-1} \\ 
&=\sum_{s\in R}
\sum_{y\in X_s }
\sum_{q=1}^{d^{\ep(r_1),\ep(r_2),\ep(y)}}
c_{r_1,r_2}^{s,q}
\;
f_{y,s}
\;
\bar\tau^{\B}_{y}\circ
m^{\ep\left(r_1\right),\ep\left(r_2\right),\ep\left(y\right)}_{q}\circ
\left(\bar\tau^{\A}_{r_1}\otimes\bar\tau^{\A}_{r_2}\right)^{-1} \\ 
&=
\sum_{y\in X }
\sum_{s\in R_y}
\sum_{q=1}^{d^{\ep(r_1),\ep(r_2),\ep(y)}}
c_{r_1,r_2}^{s,q}
\;
f_{y,s}
\;
\bar\tau^{\B}_{y}\circ
m^{\ep\left(r_1\right),\ep\left(r_2\right),\ep\left(y\right)}_{q}\circ
\left(\bar\tau^{\A}_{r_1}\otimes\bar\tau^{\A}_{r_2}\right)^{-1}. 
\end{align*}
Note that, in the 3rd equality, we replaced $\ep(s)$ by $\ep(y)$ since $y\in X_s$.

On the other hand, we have
\begin{align*}
\mu^{\B}\circ (\varphi_{r_1}\otimes \varphi_{r_2})
&=
\sum_{x_1\in X_{r_1}}
\sum_{x_2\in X_{r_2}}
f_{x_1,r_1}\, f_{x_2,r_2}\, 
\mu^{\B}_{x_1,x_2}
\left(
\bar\tau^{\B}_{x_1}\circ(\bar\tau^{\A}_{r_1})^{-1}
\otimes 
\bar\tau^{\B}_{x_2}\circ(\bar\tau^{\A}_{r_2})^{-1}
\right) \\
&=
\sum_{x_1\in X_{r_1}}
\sum_{x_2\in X_{r_2}}
\sum_{y\in X}
\sum_{q=1}^{d^{\ep(x_1),\ep(x_2),\ep(y)}}
f_{x_1,r_1}\, f_{x_2,r_2}\, 
d_{x_1,x_2}^{y,q}\; \\
&\hspace{3cm}
\bar\tau^{\B}_{y}\circ
m^{\ep\left(x_1\right),\ep\left(x_2\right),\ep\left(y\right)}_{q}\circ
\left(\bar\tau^{\B}_{x_1}\otimes\bar\tau^{\B}_{x_2}\right)^{-1}
\left(
\bar\tau^{\B}_{x_1}\circ(\bar\tau^{\A}_{r_1})^{-1}
\otimes 
\bar\tau^{\B}_{x_2}\circ(\bar\tau^{\A}_{r_2})^{-1}
\right)  \\ 
&=
\sum_{x_1\in X_{r_1}}
\sum_{x_2\in X_{r_2}}
\sum_{y\in X}
\sum_{q=1}^{d^{\ep(x_1),\ep(x_2),\ep(y)}}
f_{x_1,r_1}\, f_{x_2,r_2}\, 
d_{x_1,x_2}^{y,q}\; \
\bar\tau^{\B}_{y}\circ
m^{\ep\left(x_1\right),\ep\left(x_2\right),\ep\left(y\right)}_{q}\circ
\left(\bar\tau^{\A}_{r_1}\otimes\bar\tau^{\A}_{r_2}
\right)^{-1}\\ 
&=
\sum_{x_1\in X_{r_1}}
\sum_{x_2\in X_{r_2}}
\sum_{y\in X}
\sum_{q=1}^{d^{\ep(r_1),\ep(r_2),\ep(y)}}
f_{x_1,r_1}\, f_{x_2,r_2}\, 
d_{x_1,x_2}^{y,q}\; \
\bar\tau^{\B}_{y}\circ
m^{\ep(r_1),\ep(r_2),\ep(y)}_{q}\circ
\left(\bar\tau^{\A}_{r_1}\otimes\bar\tau^{\A}_{r_2}
\right)^{-1}.
\end{align*}
Again, in the last equality, we replaced $\ep(x_i)$ by $\ep(r_i)$. 
 
 Since the $G$-module morphisms
 \[
 \bar\tau^{\B}_{y}\circ
m^{\ep(r_1),\ep(r_2),\ep(y)}_{q}\circ
\left(\bar\tau^{\A}_{r_1}\otimes\bar\tau^{\A}_{r_2}
\right)^{-1}:
\A_{r_1}\otimes\A_{r_2}\to \B_{y},
\quad
q=1,\dots,d^{\ep(r_1),\ep(r_2),\ep(y)},
\]
are linearly independent, it follows that 
$\varphi$ is a $G$-algebra morphism if and only if 
the scalars 
$f_{x,r}\in\K$ satisfy
\[
\sum_{s\in R_y}
c_{r_1,r_2}^{s,q}
\,
f_{y,s}
=
\sum_{x_1\in X_{r_1}}
\sum_{x_2\in X_{r_2}}
d_{x_1,x_2}^{y,q}\, 
f_{x_1,r_1}\, f_{x_2,r_2}
\]
for all $r_1,r_2\in R$, $y\in X$ and 
$q=1,\dots,d^{\ep(r_1),\ep(r_2),\ep(y)}$.
\end{proof}

We close this section pointing out that 
\eqref{eq.fxr_condition} looks like as the equation defining 
morphisms of plain algebras. 
This allows us to define the plain algebras associated to 
a $G$-algebra.

\begin{Def}\label{def.plain_algebra} 
Let $G$ be a group satisfying the conditions required in \S\ref{subsec.Lableling} and assume that $\Rep(G)$ is labeled.
\begin{enumerate}
\item A \emph{choice in the labeling} of $\Rep(G)$ is 
a mapping 
$Q:(i_1,i_2,j)\mapsto Q(i_1,i_2,j)$, defined for all $i_1,i_2,j\in \hat I$ 
such that $d^{i_1,i_2,j}\ge 1$, 
where $Q(i_1,i_2,j)$ is an integer between  1 and 
$d^{i_1,i_2,j}$. 
In fact, we should think that $Q$ represents a choice of a distinguished element in $\bas^{i_1,i_2,j}$ for each 
$i_1,i_2,j\in \hat I$.
If $d^{i_1,i_2,j}\le 1$ for all $i_1,i_2,j\in \hat I$, then there is a unique possible choice in the labeling. 
\item 
Let $\A$ be a $G$-algebra with 
$G$-module decomposition $\A=\oplus_{r\in R}\A_r$ and 
$G$-table 
(see Remark \ref{rmk.abuse:language})
\begin{equation*}
    \A_{r_1}\cdot \A_{r_2}=\sum_{s\in R}\sum_{q=1}^{d^{\ep(r_1),\ep(r_2),\ep(s)}}c_{r_1,r_2}^{s,q}\;\A_{s,q}.
\end{equation*}
Let 
$P(\A)$ be the free $\K$-vector space on the set
$\{e^{\A}_r:r\in R\}$.
Given a choice $Q$ in the labeling of $\Rep(G)$ 
the \emph{plain algebra}  $P(\A,Q)$ associated to
the $G$-algebra $\A$ and 
$Q$ 
is the vector space $P(\A)$
equipped with the algebra structure 
given by the multiplication table
\begin{equation*}
    e^{\A}_{r_1}\cdot e^{\A}_{r_2}=\sum_{s\in R}c_{r_1,r_2}^{s,Q(\ep(r_1),\ep(r_2),\ep(s))}\;e^{\A}_{s}.
\end{equation*}
\end{enumerate}
\end{Def}
In this context, we have
the following functoriality property.
Let 
$\A=\oplus_{r\in R}\A_r$ and 
$\B=\oplus_{x\in X}\B_x$ be  two $G$-algebras, 
and let 
$\varphi:\A\to\B$ be a $G$-algebra morphism given by 
\begin{equation*}
\varphi_r=
\varphi|_{\A_r}=
\sum_{x\in X_r}
f_{x,r}\, \bar\tau^{\B}_{x}\circ(\bar\tau^{\A}_{r})^{-1}.
\end{equation*}
Set $f_{x,r}=0$ if $x\not\in X_r$ and let $P(\varphi)$
be the linear map given by 
\begin{align*}
P(\varphi):P(\A)&\to P(\B) \\
e^{\A}_r&\mapsto \sum_{x\in X} f_{x,r}e^{\B}_x
\end{align*}
The following corollary is a straightforward consequence of  Theorem \ref{thm.G-morphism}.
\begin{Cor}
Let $\varphi:\A\to\B$ be a $G$-module morphism.
Then $\varphi$ is a $G$-algebra morphism if and only if 
$P(\varphi):P(\A,Q)\to P(\B,Q)$ is a plain algebra morphism for any  choice $Q$ in the labeling of $\Rep(G)$.
\end{Cor}

\begin{Ej}  \label{ex.C[x]}
When $G=SL(2,\K)$, since $d^{i_1,i_2,j}\le 1$ for all 
$i_1,i_2,j\in I$, it follows that there is a unique choice $Q$
in any (partial) labeling of $\Rep(G)$, and thus a unique 
plain algebra associated to a given $G$-algebra.
In particular, the plain algebra associated
associated to $\A=\K[x,y]$ is $P(\A)=\K[x]$ 
(see Example \ref{ex.A=K[x,y]}). 
\end{Ej}

\section{The Poisson algebra associated to the cohomology of Heisenberg Lie algebra}\label{sec.Poisson_Struct}
In this section we will describe the Poisson algebra structure of 
the even cohomology $H_E^{\bullet,\bullet}(\h)$ associated to 
the cotangent bundle of 
the 3-dimensional Heisenberg Lie group with Lie algebra $\h$.
In particular we will show $H_E^{\bullet,\bullet}(\h)$ is 
the member corresponding to $n=3$ of a 
family of Poisson algebras whose associated 
Lie algebra structure is the semidirect product
$\Gl(n,\K)\ltimes\Gl(n,\K)_{ab}$, with the first factor acting on the second (abelian) one by the adjoint representation. 
Computing the $\SL(n,\K)$-table of the associative and Lie products of $H_E^{\bullet,\bullet}(\h)$ is the main tool of our proof.

\subsection{The Poisson super-algebra associated to a Lie algebra}
\label{subsec.PSA}
Given a finite dimensional Lie algebra $\g$, let
\[
C^{p,q}=C^p(\g,\bwedge^{q}\g)=\bwedge^{p}\g^{*}\otimes\bwedge^{q}\g.
\]
We know that the space 
$C^{\bullet,\bullet}=\bigoplus_{p,q}C^{p,q}$ is a 
super-commutative\footnote{Note that 
$
\left(\phi\otimes v\right)\vee\left(\psi\otimes w\right)
=
(-1)^{(p_1+q_1)(p_2+q_2)}
\left(\psi\otimes w\right)\vee
\left(\phi\otimes v\right)$} 
graded algebra 
with associative product 
\[
\vee:C^{p_1,q_1}\times C^{p_2,q_2}\rightarrow C^{p_1+p_2,q_1+q_2}
\]
 given by 
\begin{equation}\label{eq: Prod-koszul-resum}
\left(\phi\otimes v\right)\vee\left(\psi\otimes w\right)
=\left(-1\right)^{p_2q_1}\phi\wedge\psi\otimes v\wedge w.
\end{equation}
(we will refer to the factor $(-1)^{p_2q_1}$ in the above product as the Koszul sign). Moreover, 
$C^{\bullet,\bullet}$ becomes a Poisson super-algebra where the Lie super-bracket 
\[
\left\{- ,-\right\}: C^{p_1,q_1}\times C^{p_2,q_2}\rightarrow C^{p_1+p_2-1,q_1+q_2-1}
\]
is determined by (\cite{KS96})
\begin{enumerate}[label=\roman*]
\item $\left\{\g,\g\right\} =\left\{\g^{*},\g^{*}\right\} \,=0$; \label{ax 1} 
\item $\left\{a,b\right\} =-\left(-1\right)^{(p_1+q_1)(p_2+q_2)}\left\{ b,a\right\} $, 
for  $a\in C^{p_1,q_1}$, $b\in C^{p_2,q_2}$;
\label{ax 2} 
\item $\left\{ a\vee b,c\right\} =a\vee\left\{ b,c\right\} +\left(-1\right)^{(p_1+q_1)(p_2+q_2)}b\vee\left\{ a,c\right\} $, 
for  $a\in C^{p_1,q_1}$, $b\in C^{p_2,q_2}$, $c\in C^{p_3,q_3}$;
\label{ax 3} 
\item $\left\{ \phi,v\right\} =\left\{ v,\phi\right\}=\phi(v) $, 
for  $\phi\in C^{1,0}$, $v\in C^{0,1}$,
\label{ax 4} 
\end{enumerate}
Property \eqref{ax 3} is known as the {\it Poisson identity}. Since $\{-,-\}$ is a Lie super-bracket, it satisfies the super-Jacobi identity:
\begin{align}\label{eq: super-Jacobi id}
    \{a,\{b,c\}\} = \{\{a,b\},c\} + (-1)^{(p_1+q_1)(p_2+q_2)}\{b,\{a,c\}\}
\end{align}
for all $a\in C^{p_1,q_1}$, $b\in C^{p_2,q_2}$, $c\in C^{p_3,q_3}$. By the super-commutativity property \eqref{ax 2}, the item \eqref{ax 3} is equivalent to saying that $\{c,-\}$ is 
a super-derivation of the associative product $\vee$. That is
\[
\{c,a\vee b\} = \{c,a\}\vee b + (-1)^{(p_1+q_1)(p_3+q_3)}a\vee\{c,b\}
\]
for all $a\in C^{p_1,q_1}$, $b\in C^{p_2,q_2}$, $c\in C^{p_3,q_3}$. The Poisson structure on $C^{\bullet,\bullet}$ only depends on the vector space structure of $\g$ (\cite[Sec. 4]{KS96})

If $\mu\in C^{2,1}$ denotes the Lie bracket on the Lie algebra 
$\g$, then the Chevalley-Eilenberg differential $d$ on 
$C^{\bullet,\bullet}$, that computes 
$H^p(\g,\bwedge^{q}\g)$, 
the cohomology of $\g$ with coefficients on $\bwedge^{q}\g$, is 
\[
d=\{\mu,-\}:C^{p,q}\to C^{p+1,q}.
\]
Since $\{\mu,-\}$ super-derives the product $\vee$ (by \eqref{ax 3}) and the Poisson bracket $\{-,-\}$ (by \eqref{eq: super-Jacobi id}), it is well know that the cohomology space
\[
H^{\bullet,\bullet}(\g)=\bigoplus_{p,q}H^{p}\left(\mathfrak{g},\bwedge^{q}\mathfrak{g}\right)
\]
inherits, from $C^{\bullet,\bullet}$,
 the Poisson superalgebra structure and its super-commutative product $\vee$ is 
usually called the cup product.
In particular, the even subspace 
$H_E^{\bullet,\bullet}(\g)=\bigoplus_{p+q\text{ even}}H^{p}\left(\mathfrak{g},\bwedge^{q}\mathfrak{g}\right)$
is a Poisson algebra.

\subsection{The Heisenberg Lie algebra}
Let $\h$ be the Heisenberg Lie algebra of dimension 3
over a field $\K$ of characteristic zero. 
Let 
\[
B=\left\{ x_{1},x_{-1},h_{0}\right\} 
\]
be a basis of 
 $\h$ in which the only non-zero bracket is 
 $\left[x_{1},x_{-1}\right]=h_{0}$. 
We know that $G=\SL(2,\K)$ acts by automorphisms on
$\h$ and thus $\h$ is a $G$-algebra. 
Indeed, if $\{E,H,F\}$ denotes the standard $\mathfrak{s}$-triple of $\Sl(2,\K)$, then the subscripts in the basis 
$B$ indicate the $H$-weight.
The $G$-module decomposition of $\h$ is 
\[
\h=\h_0\oplus\h_1
\]
with $\h_0$ and $\h_1$ isomorphic, respectively,
to the trivial and the canonical representations of $\SL(2,\K)$.
Here we will use the partial labeling of $\Rep(\SL(2,\K))$ 
given in 
Example \ref{ex.6J-labeling}. 
In this context, if we set 
\[
\tau_0(1_{\K_0})=h_0\qquad\text{ and } \qquad
\tau_1(1,0)=x_1,\quad \tau_1(0,1)=x_{-1},
\]
then the corresponding $G$-table is 
\begin{center}
\begin{tabular}{|c||c|c|}
\hline 
 & $\h_{0}$ & $\h_1$ \\
\hline 
\hline 
$\h_0$ & $0$ &  $0$ \rule[-2mm]{0mm}{7mm} \\
\hline 
$\h_1$ & $0$ & $ \h_0$ 
\rule[-2mm]{0mm}{7mm} \\
\hline 
\end{tabular}\,.
\end{center}

Let  $B^*=\left\{ x^{-1},x^{1},h^{0}\right\} $ be the basis of $\h^*$ that is dual to $B$. 
In particular $x^{-j}(x_j)=1$, $j=\pm1$, and the superscripts in $B^*$ also indicate the $H$-weight. We additionally point out that the action of $E$ and $F$ in this dual basis is given by 
\begin{align*}
    E.x^{1}&=0,& F.x^{1}&=-x^{-1},\\
    E.x^{-1}&=-x^{1},& F.x^{-1}&=0, 
\end{align*}
and $ E.h^{0}= F.h^{0}=0$.

In \cite{CAGLIERO2004276}, the $\text{Sp}(m,\K)$-module structure of the cohomology $H^{\bullet,\bullet}(\mathfrak{H}_m)$ was given for the Heisenberg Lie algebra $\mathfrak{H}_m$ of dimension $2m+1$. 
Here, we recover from  it
the $\SL(2,\K)$-module structure of 
 $H_E^{\bullet,\bullet}(\h)$ (it corresponds to the case $m=1$ in \cite{CAGLIERO2004276}). 
 The following table 
 provides 
 the  highest weight vectors with respect to Borel subalgebra 
 of upper triangular matrices in $\Sl(2,\K)$. 
Note that we have omitted the wedge symbols $\wedge$ and that the subindices in $H^{p,q}_k$ 
indicate the highest weight of the corresponding $\SL(2,\K)$-submodule.

\begin{table}[H]
\begin{center}
\begin{tabular}{|l|l|l|}
\hline 
Structure of $H^{i,j}$ & h.w.v. of the 1st submodule & h.w.v. of the 2nd submodule \rule[-2mm]{0mm}{7mm} \\[1mm]
\hline 
\hline 
$H^{0,0}= H_{0}^{0,0}$ & $1\otimes 1$ &  \rule[-2mm]{0mm}{7mm} \\[1mm]
\hline 
$H^{2,0}= H_{1}^{2,0}$ & $x^{1}h^{0}\otimes 1$ &  \rule[-2mm]{0mm}{7mm} \\[1mm]
\hline 
$H^{1,1}= H_{0}^{1,1}\oplus H_{2}^{1,1}$ & $x^{1}\otimes x_{-1}+x^{-1}\otimes x_{1}+2h^{0}\otimes h_{0}$ & $x^{1}\otimes x_{1}$   \rule[-2mm]{0mm}{7mm} \\[1mm]
\hline $H^{3,1}= H_{1}^{3,1}$ & $x^{-1}x^{1}h^{0}\otimes x_{1}$ &   \rule[-2mm]{0mm}{7mm} \\[1mm]
\hline $H^{0,2}= H_{1}^{0,2}$ & $1\otimes x_{1}h_{0}$ &  \rule[-2mm]{0mm}{7mm} \\[1mm]
\hline $H^{2,2}= H_{0}^{2,2}\oplus H_{2}^{2,2}$ & 
$x^{-1}x^{1}\otimes x_{1}x_{-1}$ & $x^{1}h^{0}\otimes x_{1}h_{0}$  \rule[-2mm]{0mm}{7mm} \\[1mm]
\hline $H^{1,3}= H_{1}^{1,3}$ & $x^{1}\otimes x_{1}x_{-1}h_{0}$ &   \rule[-2mm]{0mm}{7mm} \\[1mm]
\hline $H^{3,3}= H_{0}^{3,3}$ & $x^{-1}x^{1}h^{0}\otimes x_{1}x_{-1}h_{0}$ &  \rule[-2mm]{0mm}{7mm} \\[1mm]
\hline 
\end{tabular}
\caption{Highest weight vectors in $H_E^{\bullet,\bullet}(\h)$}
\label{table:hwv in H}
\end{center}
\end{table}
We now give the $G$-tables of the cup product and the Poisson bracket of  $H_E^{\bullet,\bullet}(\h)$.   

\begin{thm}\label{thm.CupProduc}
Let  the $\SL(2,\K)$-morphisms $\tau$ be 
so that they send the highest weight vector $\bar v$ chosen in 
Example \ref{ex.6J-labeling} to the highest weight vector,
in the corresponding cohomology space $H_k^{p,q}(\h)$, 
given in Table \ref{table:hwv in H}. 
Then, the $G$-table for the cup product on $H_{E}^{\bullet,\bullet}$ is:
{\small{}
\begin{equation*}
\global\long\def\arraystretch{1.7}%
\begin{array}{c||ccccc|ccccc}
\vee & H_{0}^{0,0} & H_{0}^{1,1} & H_{2}^{1,1} & H_{1}^{2,0} & H_{1}^{0,2} & H_{0}^{2,2} & H_{2}^{2,2} & H_{1}^{3,1} & H_{1}^{1,3} & H_{0}^{3,3}\\
\hline\hline H_{0}^{0,0} & H_{0}^{0,0} & H_{0}^{1,1} & H_{2}^{1,1} & H_{1}^{2,0} & H_{1}^{0,2} & H_{0}^{2,2} & H_{2}^{2,2} & H_{1}^{3,1} & H_{1}^{1,3} & H_{0}^{3,3}\\
H_{0}^{1,1} & H_{0}^{1,1} & -6H_{0}^{2,2} & -2H_{2}^{2,2} & H_{1}^{3,1} & -H_{1}^{1,3} & 2H_{0}^{3,3}\\
H_{2}^{1,1} & H_{2}^{1,1} & -2H_{2}^{2,2} & H_{0}^{2,2} & H_{1}^{3,1} & H_{1}^{1,3} &  & -H_{0}^{3,3}\\
H_{1}^{2,0} & H_{1}^{2,0} & H_{1}^{3,1} & H_{1}^{3,1} &  & \frac{1}{2}H_{0}^{2,2}-\frac{1}{2}H_{2}^{2,2} &  &  &  & -H_{0}^{3,3}\\
H_{1}^{0,2} & H_{1}^{0,2} & -H_{1}^{1,3} &H_{1}^{1,3} & -\frac{1}{2}H_{0}^{2,2}-\frac{1}{2}H_{2}^{2,2} &  &  &  & -H_{0}^{3,3}\\
\hline H_{0}^{2,2} & H_{0}^{2,2} & 2H_{0}^{3,3} &  &  & \\
H_{2}^{2,2} & H_{2}^{2,2} &  & -H_{0}^{3,3} &  & \\
H_{1}^{3,1} & H_{1}^{3,1} &  &  &  & H_{0}^{3,3}\\
H_{1}^{1,3} & H_{1}^{1,3} &  &  & H_{0}^{3,3} & \\
H_{0}^{3,3} & H_{0}^{3,3} &  &  &  & 
\end{array}
\end{equation*}}
\end{thm}

\begin{remark}
Before proceeding with the proof, we point out that, even though the cup product is commutative in $H_E^{\bullet,\bullet}(\h)$, there are some entries that may look as indicating the opposite. 
They are
\begin{align*}
H_{1}^{2,0}\vee H_{1}^{0,2}&=\frac{1}{2}H_{0}^{2,2} -\frac{1}{2}H_{2}^{2,2}& \text{ while }&&
H_{1}^{0,2}\vee H_{1}^{2,0}&=-\frac{1}{2}H_{0}^{2,2} -\frac{1}{2}H_{2}^{2,2}   \\
H_{1}^{2,0}\vee H_{1}^{1,3}&=- H_{0}^{3,3}& \text{ while }&&
H_{1}^{1,3}\vee H_{1}^{2,0} &= H_{0}^{3,3}\\ 
H_{1}^{0,2}\vee H_{1}^{3,1}&=- H_{0}^{3,3}& \text{ while }&&
H_{1}^{3,1}\vee H_{1}^{0,2} &= H_{0}^{3,3}.
\end{align*}
All of these are due to the fact that the $\SL(2,\K)$-morphism
$m^{1,1,0}:(\pi_{1},\K^2)\otimes (\pi_{1},\K^2)\to (\pi_{0},\K_{0})$ is skew-symmetric (see Example \ref{ex.6J-labeling}).
\end{remark}

\begin{proof}
In all the following computations it is important to recall the choice of the maps $\tau$ that we have already made. 
It is also useful to recall that $F.x_1=x_{-1}$ and $F.x^1=-x^{-1}$.

We will organize the proof following the rows of Table \ref{table:hwv in H}. Since the cup product in $H_E^{\bullet,\bullet}(\h)$
is commutative, it is enough to consider the upper right triangle 
of the table.
Also, since the zero products are easier to check than the non-zero ones (many of them are straightforward because there is no room for the potential result) we will concentrate in the non-zero entries. 

\begin{itemize}
    \item Proof for the row  corresponding to $H_0^{0,0}$. 
    Since $H_0^{0,0}$ is the module spanned by $1\otimes 1$ which is 
    the identity of the cup product, the proof is straightforward. 
    \item Proof for the row  corresponding to $H_0^{1,1}$. 
    Here, the most subtle computation corresponds to the entry is 
    $H_0^{1,1}\vee H_0^{1,1} = -6\,H_0^{2,2}$. To prove it, we have
    \[
       \big(x^{1}\otimes x_{-1}+x^{-1}\otimes x_{1}+2h^{0}\otimes h_{0}\big)\vee\big(x^{1}\otimes x_{-1}+x^{-1}\otimes x_{1}+2h^{0}\otimes h_{0}\big)
        =-6\,x^{-1}x^{1}\otimes x_{1}x_{-1} +4 u    
    \]
    with 
    \[
    u= x^{-1}x^{1}\otimes x_{1}x_{-1} -x^{-1}h^{0}\otimes x_{1}h_{0}- x^{1}h^{0}\otimes x_{-1}h_{0}
    = d(h^{0}\otimes x_{1}x_{-1}).
    \]
    On other hand, since ${m}^{0,0,0}(1_{\K_0}\otimes 1_{\K_0})=1_{\K_0}$, it follows that
    \[
        \tilde{m}^{0,0,0}\big(x^{1}\otimes x_{-1}+x^{-1}\otimes x_{1}+2h^{0}\otimes h_{0}\, ,\,x^{1}\otimes x_{-1}+x^{-1}\otimes x_{1}+2h^{0}\otimes h_{0}\big) = x^{-1}x^{1}\otimes x_{1}x_{-1}
    \]
    and thus the cup product in $H_0^{1,1}\otimes H_0^{1,1}$ is 
    $-6\,\tilde{m}^{0,0,0}$.
    
    The other entries of this row:
    \begin{align*}
        H_0^{1,1}\vee H_2^{1,1} &= -2H_2^{2,2}, &
        H_0^{1,1}\vee H_1^{2,0} &= H_1^{3,1}, \\
        H_0^{1,1}\vee H_1^{0,2} &= -H_1^{1,3}, &
        H_0^{1,1}\vee H_0^{2,2} &= 2H_0^{3,3},
    \end{align*}
 are similar, but simpler. We point out that the minus sign in 
$H_0^{1,1}\vee H_2^{1,1} = -2H_2^{2,2}$ comes from the Koszul sign $(-1)^{p_2q_1}$ in \eqref{eq: Prod-koszul-resum}.

\item Proof for the row  corresponding to $H_2^{1,1}$. 
\begin{itemize}
    \item Proof of $H_2^{1,1}\vee H_1^{0,2} = H_1^{1,3}$. 
    It is clear that we must have $H_2^{1,1}\vee H_1^{0,2} = c\, H_1^{1,3}$ for some scalar $c$. This $c$ can be obtained by computing a non-zero (if any) cup product in the corresponding modules. 
    On the one hand, we have
    \begin{align*}
     (x^1\otimes x_1)\vee F. (1\otimes x_{1}h_0) & = (x^1\otimes x_1)\vee  (1\otimes x_{-1}h_0) \\ 
     &= x^1\otimes x_1x_{-1}h_0.
    \end{align*} 
    On the other hand, we know that (see Example \ref{ex.6J-labeling})
    \[
    m^{2,1,1}\left( 
    \left( \begin{smallmatrix}
        0 & 1 \\ 0 & 0 
    \end{smallmatrix} \right)\otimes F. (1,0)\right) = 
    m^{2,1,1}\left( 
    \left( \begin{smallmatrix}
        0 & 1 \\ 0 & 0 
    \end{smallmatrix} \right)\otimes (0,1)\right) = (1,0)
    \]
    and thus $c=1$.

    \item Proof of $H_2^{1,1}\vee H_1^{2,0} = H_1^{3,1}$. The argument here is almost the same as in the previous case, one detail to be aware of is that $F. (x^{1}h_0 \otimes 1)=-(x^{-1}h_0 \otimes 1)$.
    
      \item Proof of $H_2^{1,1}\vee H_2^{1,1} =  H_0^{2,2}$. 
    At first glance, we must have 
    \[
    H_2^{1,1}\vee H_2^{1,1} = c_0\, H_0^{2,2} + c_2\, H_2^{2,2}
    \] for some scalars $c_0,c_2$. 
    But, since the cup product is commutative and 
    $m^{2,2,2}$ is skew-symmetric ($m^{2,2,2}(A\otimes B)=-m^{2,2,2}(B\otimes A)$, see Example \ref{ex.6J-labeling}), we conclude that $c_2=0$. 
    To obtain $c_0$ we need to compute a non-zero (if any) cup product in  $H_2^{1,1}$. 
    On the one hand, we have
    \begin{align*}
     (x^1\otimes x_1)\vee F^2. (x^1\otimes x_1) 
     & = -2(x^1\otimes x_1)\vee (x^{-1}\otimes x_{-1}) \\ 
     & =  2(x^1x^{-1}\otimes x_1x_{-1}) \qquad \text{(note the Koszul sign)} \\
    & =  -2(x^{-1}x^1\otimes x_1x_{-1}).
    \end{align*} 
    On the other hand, we know that (see Example \ref{ex.6J-labeling})
    \[
    m^{2,2,0}\left( 
    \left( \begin{smallmatrix}
        0 & 1 \\ 0 & 0 
    \end{smallmatrix} \right)\otimes 
    F^2.
    \left(\begin{smallmatrix}
        0 & 1 \\ 0 & 0 
    \end{smallmatrix} \right)
    \right) = 
    m^{2,2,0}\left( 
    \left( \begin{smallmatrix}
        0 & 1 \\ 0 & 0 
    \end{smallmatrix} \right)\otimes 
   \left( \begin{smallmatrix}
        0 & 0 \\ -2 & 0 
    \end{smallmatrix} \right)
    \right) 
    = \text{tr} \left(\begin{smallmatrix}
        -2 & 0 \\ 0 & 0 
    \end{smallmatrix} \right)
    =-2
    \]
    and thus $c_0=1$.

      \item Proof of $H_2^{1,1}\vee H_2^{2,2} =  -H_0^{3,3}$. 
      It is almost the same as the previous case, the only difference 
      is that, here, the Koszul sign is positive instead of negative.
\end{itemize}

    \item  Proof for the row  corresponding to $H_1^{2,0}$. 
 \begin{itemize}
   \item Proof of  $H_1^{2,0}\vee H_1^{0,2} = \frac{1}{2}H_0^{2,2}-\frac{1}{2}H_2^{2,2}$. 
   We know that $H_1^{2,0}\vee H_1^{0,2} = c_0H_0^{2,2}+c_2H_2^{2,2}$
   for certain scalars $c_0,c_2$. 
   We have, on the one hand,
   \begin{align}
   \label{eq.aux11}
   (x^1h^0\otimes 1)\vee (1\otimes x_1h_0) & = x^1h^0\otimes x_1h_0 \\
   (x^1h^0\otimes 1)\vee F.(1\otimes x_1h_0) & = x^1h^0\otimes x_{-1}h_0\notag  \\
   & = \tfrac12\big(
   F.( x^1h^0\otimes x_1h_0) - d(h^{0}\otimes x_{1}x_{-1}) 
   + x^{-1}x^{1}\otimes x_{1}x_{-1}\big)\notag \\
    \label{eq.aux12}
  & = \tfrac12
   F.( x^1h^0\otimes x_1h_0) 
   + \tfrac12(x^{-1}x^{1}\otimes x_{1}x_{-1}). 
   \end{align}
   Note that 
   \[
   d(h^{0}\otimes x_{1}x_{-1})= 
    x^{-1}x^{1}\otimes x_{1}x_{-1}
   -x^{-1}h^{0}\otimes x_{1}h_{0} 
   -x^{1}h^{0}\otimes x_{-1}h_{0}.
   \]
   On the other hand  (see Example \ref{ex.6J-labeling}),
   \begin{align}
   \label{eq.aux21}
   m^{1,1,0}\big((1,0)\otimes (1,0)\big)&=0, &
    m^{1,1,2}\big((1,0)\otimes (1,0)\big)&=
   \left(\begin{smallmatrix} 0&-2\\0&0\end{smallmatrix}\right),  \\
   \label{eq.aux22} 
   m^{1,1,0}\big((1,0)\otimes F.(1,0)\big)&=1_{\K_0} &
   m^{1,1,2}\big((1,0)\otimes F.(1,0)\big)&=
   \left(\begin{smallmatrix} 1&0\\0&-1\end{smallmatrix}\right)=
   -F. \left(\begin{smallmatrix} 0&1\\0&0\end{smallmatrix}\right).
   \end{align}
From \eqref{eq.aux11} and \eqref{eq.aux21}, we obtain $c_2=-\frac12$.
Next, from \eqref{eq.aux12} and \eqref{eq.aux22}, we obtain $c_0=\frac12$.

    \item Proof of $H_1^{2,0}\vee H_1^{1,3} = -H_0^{3,3}$.
    It is clear that we must have 
    $H_1^{2,0}\vee H_1^{1,3} = c H_0^{3,3}$ for some scalar $c$. 
    To obtain $c$ we have, on the one hand, we have
    \begin{align*}
     (x^1h^0\otimes 1)\vee F. (x^1\otimes x_1x_{-1}h_0) & 
     = -(x^1h^0\otimes 1)\vee  ( x^{-1}\otimes x_1x_{-1}h_0) \\ 
     &= - (x^{-1}x^1h^0\otimes x_1x_{-1}h_0).
    \end{align*} 
    On the other hand, we know that (see Example \ref{ex.6J-labeling})
    \[
    m^{1,1,0}\big( (1,0)\otimes F.(1,0) \big) =1_{\K_0}
    \]
    and thus $c=-1$. 
\end{itemize}
\item  Proof for the row  corresponding to $H_1^{0,2}$. 
The only cup product that has not already been computed for this row is
$H_1^{0,2}\vee H_1^{3,1} = -H_0^{3,3}$, and the proof of it is very similarly to that of the entry $H_1^{2,0}\vee H_1^{1,3} = -H_0^{3,3}$. 
\end{itemize}
This completes the proof. 
\end{proof}

\begin{thm}\label{thm.PoissonBracket}
With the same choices made for Theorem \ref{thm.CupProduc}, the $G$-table of the Poisson bracket on $H_{E}^{\bullet,\bullet}$ is:
{\scriptsize{}
\begin{equation*}
\global\long\def\arraystretch{1.8}%
\begin{array}{c||ccccc|ccccc}
\left\{ -,-\right\}  & H_{0}^{0,0} & H_{0}^{1,1} & H_{2}^{1,1} & H_{1}^{2,0} & H_{1}^{0,2} & H_{0}^{2,2} & H_{2}^{2,2} & H_{1}^{3,1} & H_{1}^{1,3} & H_{0}^{3,3}\\
\hline\hline H_{0}^{0,0} &  &  &  &  & \\
H_{0}^{1,1} &  &  &  & 3H_{1}^{2,0} & -3H_{1}^{0,2} &  &  & 3H_{1}^{3,1} & -3H_{1}^{1,3}\\
H_{2}^{1,1} &  &  & -H_{2}^{1,1} & -H_{1}^{2,0} & -H_{1}^{0,2} &  & -H_{2}^{2,2} & -H_{1}^{3,1} & -H_{1}^{1,3}\\
H_{1}^{2,0} &  & -3H_{1}^{2,0} & H_{1}^{2,0} &  & \frac{1}{2}H_{2}^{1,1}-\frac{1}{2}H_{0}^{1,1} & H_{1}^{3,1} & H_{1}^{3,1} &  & -\frac{1}{2}H_{2}^{2,2}-\frac{3}{2}H_{0}^{2,2}\\
H_{1}^{0,2} &  & 3H_{1}^{0,2} & H_{1}^{0,2} & -\frac{1}{2}H_{2}^{1,1}-\frac{1}{2}H_{0}^{1,1} &  & H_{1}^{1,3} & -H_{1}^{1,3} &  -\frac{1}{2}H_{2}^{2,2}+\frac{3}{2}H_{0}^{2,2}\\
\hline H_{0}^{2,2} &  &  &  & -H_{1}^{3,1} & -H_{1}^{1,3}\\
H_{2}^{2,2} &  &  & -H_{2}^{2,2} & -H_{1}^{3,1} & H_{1}^{1,3}\\
H_{1}^{3,1} &  & -3H_{1}^{3,1} & H_{1}^{3,1} &  & +\frac{1}{2}H_{2}^{2,2}+\frac{3}{2}H_{0}^{2,2}\\
H_{1}^{1,3} &  & 3H_{1}^{1,3} & H_{1}^{1,3} & \frac{1}{2}H_{2}^{2,2}-\frac{3}{2}H_{0}^{2,2} & \\
H_{0}^{3,3} &  &  &  &  & 
\end{array}
\end{equation*}}
\end{thm}

\begin{proof}
As in the previous proof, we will concentrate in the non-zero entries of the table. 
Along the proof we will frequently use of Properties \eqref{ax 1}-\eqref{ax 4},  \S\ref{subsec.PSA}, for  the even part of 
the Poisson bracket $\{-,-\}:H_E^{p_1,q_1}\otimes H_E^{p_2,q_2}\rightarrow H_E^{p_1+p_2-1,q_1+q_2-1}$.
\begin{itemize}
\item Proof for the row corresponding to $H_0^{1,1}$. Let 
\[
v=x^1\otimes x_{-1} + x^{-1}\otimes x_1 + 2h^0\otimes h_0\in H_0^{1,1}.
\]
\begin{itemize}
    \item Proof of $\{H_0^{1,1}, H_1^{2,0}\} = 3H_1^{2,0}$. 
    We have 
    \[
    \{v,x^1h^0\otimes 1\} = 3\,(x^1h^0\otimes 1).
    \] 
    Indeed, for instance
    \begin{align*}
     \{x^1\otimes x_{-1},x^1h^0\otimes 1\} 
     & = (x^1\otimes 1)\vee \{1\otimes x_{-1}, x^{1}h^0\otimes 1\} - (1\otimes x_{-1})\vee \{ x^1\otimes 1, x^{1}h^0\otimes 1\} \\
     & = (x^1\otimes 1)\vee \{1\otimes x_{-1}, x^{1}h^0\otimes 1\} \\
     & = (x^1\otimes 1)\vee (h^0\otimes 1)  \\
     &=(x^1h^0\otimes 1)
        \end{align*}
    since $\{\g^*,\g^*\} = 0$ and  $\left\{ \phi,v\right\} =\left\{ v,\phi\right\}=\phi(v) $, 
for  $\phi\in C^{1,0}$, $v\in C^{0,1}$.
Similarly we obtain
    \begin{align*}
\{x^{-1}\otimes x_{1},x^1h^0\otimes 1\} &=0 \\
\{h^{0}\otimes h_{0},x^1h^0\otimes 1\} &=(x^1h^0\otimes 1),
        \end{align*}
and thus    
$\{x^1\otimes x_{-1} + x^{-1}\otimes x_1 + 2h^0\otimes h_0,x^1h^0\otimes 1\} = 3\,(x^1h^0\otimes 1)$ as claimed.

\medskip
\item Proof of $\{H_0^{1,1}, H_1^{3,1}\} = 3H_1^{3,1}$. 
First, note that 
\[
 \{(x^{-1}\otimes x_1),(x^{-1}\otimes x_1)\}=
 \{(x^{1}\otimes x_{-1}),(x^{-1}\otimes x_1)\}=
 \{(h^{0}\otimes h_0),(x^{-1}\otimes x_1)\}=0
\]
and hence $\{v,(x^{-1}\otimes x_1)\}=0$.
Therefore
    \begin{align*}
     \{v,x^{-1}x^1h^0\otimes x_1\} 
     & = \{v,(x^{-1}\otimes x_1)\vee (x^1h^0\otimes 1)\} \\
    & = \{v,(x^{-1}\otimes x_1)\}\vee (x^1h^0\otimes 1)
         +(x^{-1}\otimes x_1)\vee \{v, x^{1}h^0\otimes 1\} \\
     & = (x^{-1}\otimes x_1)\vee \{v, x^{1}h^0\otimes 1\} \\
     & = (x^{-1}\otimes x_1)\vee 3\, (x^{1}h^0\otimes 1)\qquad\text{ by the previous computation} \\
     &=3\,(x^{-1}x^1h^0\otimes x_1).
        \end{align*}    
\medskip
    \item Proof of $\{H_0^{1,1}, H_1^{0,2}\} = -3H_1^{0,2}$
    and $\{H_0^{1,1}, H_1^{1,3}\} = -3H_1^{1,3}$. These two cases are similar to the previous two. Here
    \[
    \{v,1\otimes x_1h_0\} = -3\,(1\otimes x_1h_0).
    \]
\end{itemize}
    
\medskip
    \item Proof for the row corresponding to $H_2^{1,1}$.
    \begin{itemize}
        \item Proof of $\{H_2^{1,1}, H_2^{1,1}\} = -H_2^{1,1}$. 
    We know that $\{H_2^{1,1}, H_2^{1,1}\} = c_0H_0^{1,1}+c_2H_2^{1,1}$, but $c_0$ must be zero since the bracket is skew-symmetric but 
    $m^{2,2,0}$ is symmetric. That $c_2=-1$ now follows from
    \[
    \{x^1\otimes x_1,F.(x^1\otimes x_{1})\} = 
    \{x^1\otimes x_1,x^1\otimes x_{-1}-x^{-1}\otimes x_1\} = 
    -2(x^1\otimes x_1)
    \]
    and 
    $m^{2,2,2}\left( 
    \left( \begin{smallmatrix}
        0 & 1 \\ 0 & 0 
    \end{smallmatrix} \right)\otimes 
   \left( \begin{smallmatrix}
        -1 & 0 \\ 0 & 1 
    \end{smallmatrix} \right)
    \right) 
    = 2    \left( \begin{smallmatrix}
        0 & 1 \\ 0 & 0 
    \end{smallmatrix} \right)
    $.
    Similarly it is obtained that  $\{H_2^{1,1}, H_2^{2,2}\} = -H_2^{2,2}$. 
    
\medskip
\item  Proof of $\{H_2^{1,1}, H_1^{2,0}\} = -H_1^{2,0}$. 
    This follows from 
    \begin{align*}
     \{x^1\otimes x_1,F.(x^{1}h^0\otimes 1)\} 
     &= \{x^1\otimes x_1,-x^{-1}h^0\otimes 1\} = \\ 
     & = (x^1\otimes 1)\vee \{1\otimes x_1, -x^{-1}h^0\otimes 1\} - (1\otimes x_1)\vee \{ x^1\otimes 1, -x^{-1}h^0\otimes 1\} \\
     & = -(x^1\otimes 1)\vee \{1\otimes x_1, x^{-1}h^0\otimes 1\} \\
     & = -(x^1\otimes 1)\vee (h^0\otimes 1)  \\
     &=-(x^1h^0\otimes 1)
        \end{align*}
       and     $m^{2,1,1}\left( 
    \left( \begin{smallmatrix}
        0 & 1 \\ 0 & 0 
    \end{smallmatrix} \right)\otimes F.(1,0) \right)
    = (1,0) $. 
    
    \item  Proof of  
    $\{H_2^{1,1}, H_1^{0,2}\}= -H_1^{0,2}$, 
    $\{H_2^{1,1}, H_1^{3,1}\} = -H_1^{3,1}$ and 
    $\{H_2^{1,1}, H_1^{1,3}\} = -H_1^{1,3}$. 
   All of these are similar to the proof of 
  $\{H_2^{1,1}, H_1^{2,0}\} = -H_1^{2,0}$.
  They follow from 
      \begin{align*}
     \{x^1\otimes x_1,F.(1\otimes x_{1}h_0)\} 
     &= \{x^1\otimes x_1,1 \otimes x_{-1}h_0\} \\ 
     &=-(1\otimes x_1h_0),
        \end{align*}
      \begin{align*}
     \{x^1\otimes x_1,F.(x^1\otimes x_{1}x_{-1}h_0)\} 
     &= \{x^1\otimes x_1,-x^{-1}\otimes x_{1}x_{-1}h_0\} \\ 
     &=-(x^{1}\otimes x_{1}x_{-1}h_0),
        \end{align*}
        and
      \begin{align*}
     \{x^1\otimes x_1,F.(x^{-1}x^1h^0\otimes x_{1})\} 
     &= \{x^1\otimes x_1,x^{-1}x^1h^0\otimes x_{-1}\} \\ 
     &=-(x^{-1}x^1h^0\otimes x_{1}).
        \end{align*}
     \end{itemize}
 
    \item Proof for the row corresponding to $H_1^{2,0}$.
    \begin{itemize}
    \item Proof of $\{H_1^{2,0}, H_2^{2,2}\} = H_1^{3,1}$. 
    We know that, by definition,
    \[
 m^{1,2,1}\left(F.(1,0) \otimes 
    \left( \begin{smallmatrix}
        0 & 1 \\ 0 & 0 
    \end{smallmatrix} \right) \right)
    =
m^{2,1,1}\left( 
    \left( \begin{smallmatrix}
        0 & 1 \\ 0 & 0 
    \end{smallmatrix} \right)\otimes F.(1,0) \right)
    = (1,0). 
    \]
    On the other hand,
    \begin{align*}
     \{F.(x^{1}h^0\otimes 1)&,(x^{1}h^0\otimes x_1h_0)\} \\
     & = -\{F.(x^{1}h^0\otimes 1),(x^{1}\otimes x_1)\vee (h^0\otimes h_0)\} \\
     & = -\{F.(x^{1}h^0\otimes 1),(x^{1}\otimes x_1)\}\vee (h^0\otimes h_0)      -(x^{1}\otimes x_1)\vee \{F.(x^{1}h^0\otimes 1),(h^0\otimes h_0)\} \\
     & = -(x^{1}h^0\otimes 1)\vee (h^0\otimes h_0)     
     +(x^{1}\otimes x_1)\vee \{(x^{-1}h^0\otimes 1),(h^0\otimes h_0)\} \\
     & = -(x^{1}\otimes x_1)\vee (x^{-1}h^0\otimes 1) \\
     &=(x^{-1}x^1h^0\otimes x_1),
        \end{align*}    
    here we used that
    $\{F.(x^{1}h^0\otimes 1),(x^{1}\otimes x_1)\} = (x^{1}h^0\otimes 1)$,
    it was shown when we
    proved that $\{H_2^{1,1}, H_1^{2,0}\} = -H_1^{2,0}$.
    
    \item Proof of $\{H_1^{2,0}, H_0^{2,2}\} = H_1^{3,1}$. 
    This follows from 
         \begin{align*}
     \{(x^{1}h^0\otimes 1)&,(x^{-1}x^{1}\otimes x_1x_{-1})\} \\
     & = \{(x^{1}h^0\otimes 1),(x^{1}\otimes x_1)\vee (x^{-1}\otimes x_{-1})\} \\
     & = \{(x^{1}h^0\otimes 1),(x^{1}\otimes x_1)\}\vee (x^{-1}\otimes x_{-1})
     +(x^{1}\otimes x_1)\vee \{(x^{1}h^0\otimes 1),(x^{-1}\otimes x_{-1})\} \\
     & =  (x^{1}\otimes x_1)\vee \{(x^{1}h^0\otimes 1),(x^{-1}\otimes x_{-1})\} \\
     & = -(x^{1}\otimes x_1)\vee (x^{-1}h^0\otimes 1) \\
     &=(x^{-1}x^1h^0\otimes x_1),
        \end{align*}  
           here we used that
    $\{(x^{1}h^0\otimes 1),(x^{-1}\otimes x_{-1})\} = -(x^{-1}h^0\otimes 1)$,
    this follows from the fact that 
    $\{H_2^{1,1}, H_1^{2,0}\} = -H_1^{2,0}$ and 
    $m^{2,1,1}\left( 
    \left( \begin{smallmatrix}
        0 & 0 \\ 1 & 0 
    \end{smallmatrix} \right)\otimes (1,0) \right)
    = (1,0)$.

\medskip
    \item Proof of $\{H_1^{2,0}, H_1^{0,2}\} = \frac{1}{2}H_2^{1,1} - \frac{1}{2}H_0^{1,1}$
    and 
    $\{H_1^{2,0}, H_1^{1,3}\} = -\frac{1}{2}H_2^{2,2} - \frac{1}{2}H_0^{2,2}$. 
    We will prove the first one, the second one is similar.
    We know that 
    $\{H_1^{2,0}, H_1^{0,2}\} = c_2H_2^{1,1} +c_0H_0^{1,1}$.
    Since 
    \[
    \{(x^1h^0\otimes 1),(1\otimes x_1h_0)\} = -x^1\otimes x_1
    \]
    and $m^{1,1,2}\big((1,0)\otimes (1,0)\big)
    =-2
    \left( \begin{smallmatrix}
        0 & 1 \\ 0 & 0 
    \end{smallmatrix} \right)$, it follows that $c_2=\frac12$.
    
    On the other hand, 
    \begin{align*}
    \{F.(x^{1}h^0\otimes 1),&(1\otimes x_1h_0)\} - \{(x^1h^0\otimes 1),F.(1\otimes x_{1}h_0)\} \\
    &=    
    \{-x^{-1}h^0\otimes 1,1\otimes x_1h_0\} - \{x^1h^0\otimes 1,1\otimes x_{-1}h_0\} \\
    &= x^1\otimes x_{-1} + x^{-1}\otimes x_1 + 2\,(h^0\otimes h_0),
    \end{align*}
    and, since
    \begin{align*}
    m^{1,1,0}(F.(1,0)\otimes (1,0))-m^{1,1,0}((1,0)\otimes F.(1,0))&=-2, \\
    m^{1,1,2}(F.(1,0)\otimes (1,0))-m^{1,1,0}((1,0)\otimes F.(1,0))&=\left( \begin{smallmatrix}
        0 & 0 \\ 0 & 0 
    \end{smallmatrix} \right), 
    \end{align*}
    it follows that $c_0=-\frac{1}{2}$.
    \end{itemize}

\medskip

    \item Proof for the row corresponding to $H_1^{0,2}$.
    
\begin{itemize}
     \item Proof of  $\{H_1^{0,2}, H_0^{2,2}\} = H_1^{1,3}$
     and  $\{H_1^{0,2}, H_2^{2,2}\} = H_1^{1,3}$.
     These are similar to the proofs of  
  $\{H_1^{2,0}, H_2^{2,2}\} = H_1^{3,1}$ and 
  $\{H_1^{2,0}, H_0^{2,2}\} = H_1^{3,1}$.
  
     \item  Proof of   $\{H_1^{0,2}, H_1^{3,1}\} = \frac{1}{2}H_2^{2,2} + \frac{3}{2}H_0^{2,2}$. The argument here is almost the same as in the proof of 
    $\{H_1^{2,0}, H_1^{0,2}\} = \frac{1}{2}H_2^{1,1} - \frac{1}{2}H_0^{1,1}$.
    Indeed, we know that 
    $\{H_1^{0,2}, H_1^{3,1}\} = c_2 H_2^{2,2} + c_0 H_0^{2,2}$.
    Since
    \[
        \{ 1\otimes x_1h_0, x^{-1}x^1h^0\otimes x_1\} = x^1h^0\otimes x_1h_0
    \]
    and hence $c_2=-\frac12$.
    On the other hand
    \begin{align*}
        \{ F.(1\otimes x_1h_0), (x^{-1}x^1h^0\otimes x_1)\} 
       &= \{ 1\otimes x_{-1}h_0, x^{-1}x^1h^0\otimes x_{1}\} \\
        &= -x^{-1}x^1\otimes x_1x_{-1} - x^{-1}h^0\otimes x_{1}h_0;\\
        \{ (1\otimes x_1h_0), F.(x^{-1}x^1h^0\otimes x_1)\} 
        &= \{ 1\otimes x_1h_0, x^{-1}x^1h^0\otimes x_{-1}\} \\
        &= x^{-1}x^1\otimes x_1x_{-1} + x^1h^0\otimes x_{-1}h_0
    \end{align*}
    and thus $\{ F.(1\otimes x_1h_0), (x^{-1}x^1h^0\otimes x_1)\} 
        -\{ (1\otimes x_1h_0), F.(x^{-1}x^1h^0\otimes x_1)\}$
        equals to
    \[
      -2x^{-1}x^1\otimes x_1x_{-1} - x^{-1}h^0\otimes x_{1}h_0 - 
       x^1h^0\otimes x_{-1}h_0 
     = -3x^{-1}x^1\otimes x_1x_{-1} + d(h^{0}\otimes x_{1}x_{-1}),       
    \]
    and hence $c_0=\frac{3}{2}$.
    \end{itemize}
\end{itemize}
This completes the proof. 
\end{proof}

\begin{Cor}\label{Thm.Q3-HE(h1)}
As a Lie algebra, $H_{E}^{\bullet,\bullet}(\h)$ is
isomorphic to $\Gl(3,\K)\ltimes\Gl(3,\K)_{ab}$,
where $\Gl(3,\K)_{ab}$ is viewed as an abelian Lie algebra on which $\Gl(3,\K)$ acts by the adjoint representation. 
\end{Cor}

\begin{proof} 
We can see that the left part of the $G$-table of the Poisson bracket
of $H_{E}^{\bullet,\bullet}(\h)$ 
{\scriptsize{}
\begin{equation*}
\global\long\def\arraystretch{1.8}%
\begin{array}{c||ccccc|ccccc}
\left\{ -,-\right\}  & H_{0}^{0,0} & H_{0}^{1,1} & H_{2}^{1,1} & H_{1}^{2,0} & H_{1}^{0,2} & H_{0}^{2,2} & H_{2}^{2,2} & H_{1}^{3,1} & H_{1}^{1,3} & H_{0}^{3,3}\\
\hline\hline H_{0}^{0,0} &  &  &  &  & \\
H_{0}^{1,1} &  &  &  & 3H_{1}^{2,0} & -3H_{1}^{0,2} &  &  & 3H_{1}^{3,1} & -3H_{1}^{1,3}\\
H_{2}^{1,1} &  &  & -H_{2}^{1,1} & -H_{1}^{2,0} & -H_{1}^{0,2} &  & -H_{2}^{2,2} & -H_{1}^{3,1} & -H_{1}^{1,3}\\
H_{1}^{2,0} &  & -\underline{3}H_{1}^{2,0} & H_{1}^{2,0} &  & \underline{\frac{1}{2}}H_{2}^{1,1}-\underline{\frac{1}{2}}H_{0}^{1,1} & H_{1}^{3,1} & H_{1}^{3,1} &  & -\frac{1}{2}H_{2}^{2,2}-\frac{3}{2}H_{0}^{2,2}\\
H_{1}^{0,2} &  & \underline{3}H_{1}^{0,2} & \underline{1}H_{1}^{0,2} & -\frac{1}{2}H_{2}^{1,1}-\underline{\frac{1}{2}}H_{0}^{1,1} &  & H_{1}^{1,3} & -H_{1}^{1,3} &  -\frac{1}{2}H_{2}^{2,2}+\frac{3}{2}H_{0}^{2,2}
\end{array}
\end{equation*}}
is almost the same as the right half. 
The only differences are the scalars that are underlined in the left half. 
These differences could be `repaired' by rescaling a couple of h.w.v in 
Table \ref{table:hwv in H}. 
This implies that, as a Lie algebra, 
\[
H_{E}^{\bullet,\bullet}(\h)\simeq
\g_0\ltimes\g_0
\]
where $\g_0$ is the Lie subalgebra
\[
\g_0= H_{0}^{0,0} \;\oplus\; H_{0}^{1,1} \;\oplus\; H_{2}^{1,1} \;\oplus\; H_{1}^{2,0} \;\oplus\; H_{1}^{0,2}.
 \]
 
On the other hand, the $SL(2,\K)$-table of $\mathfrak{sl}(3,\K)$
(see Example \ref{ex.A=sl(3)}) is
\begin{center}
\small
\begin{tabular}{|c||c|c|c|c|}
\hline 
 & {$V_0$} & {$V_2$} & {$V_1$} & {$V_1'$} \rule[-2mm]{0mm}{6mm} \\
\hline 
\hline 
 {$V_0$} & & & {$3\,V_1$} & {$-3\,V_1'$} \rule[-1mm]{0mm}{6mm} \\
\hline 
 {$V_2$} &  & {$V_2$} & {$V_1$} & {$V_1'$} \rule[-1mm]{0mm}{6mm} \\
\hline 
 {$V_1$} & {$-3\,V_1$} & {$-V_1$} &  & {$-\frac{1}{2}\,V_2+\frac{1}{2}\,V_0$} \rule[-2mm]{0mm}{6mm} \\
\hline 
 {$V_1'$} & {$3\,V_1$}  & {$-V_1'$} & {$\frac{1}{2}\,V_2+\frac{1}{2}\,V_0$} &  \rule[-2mm]{0mm}{6mm} \\
\hline 
\end{tabular}
\end{center}
which can also be easily rescaled to math the left half of the above $G$-table. 
This implies that $\g_0\simeq \Gl(3,\K)$.
\end{proof}
\begin{remark}
It is remarkable for us that the Lie algebra structure on 
$H_{E}^{\bullet,\bullet}(\h)$ contains a semisimple 
Lie subalgebra (in this case $\Sl(3,\K)$) 
strictly larger than the Levi factor 
of $\text{Der}(\h)$, which in this case is $\Sl(2,\K)\subset H^{1,1}$. 
\end{remark}

As a by product of this computation, we discovered that 
$H_{E}^{\bullet,\bullet}(\h)$ is a member of 
an infinite family of Poisson algebras whose underlying 
Lie algebra is the semidirect product
$\Gl(n,\K)\ltimes\Gl(n,\K)_{ab}$.
As  $\Gl(n,\K)$-module, 
$\Gl(n,\K)\ltimes\Gl(n,\K)_{ab}$ is a sum of two copies 
of the trivial module and two copies of the adjoint module. 
We choose these copies inside the factors of the semidirect product and thus, we denote these copies as 
\begin{equation}\label{eq.decomp_semidirect}
\Gl(n,\K)\ltimes\Gl(n,\K)_{ab}\simeq 
\big((I_n)_0\oplus \Sl(n,\K)_0\big)
\;\oplus\; \big((I_n)_{ab}\oplus \Sl(n,\K)_{ab}\big).
\end{equation}
Accordingly, we describe, in the following theorem, 
the elements in 
$\Gl(n,\K)\ltimes\Gl(n,\K)_{ab}$ as pairs
$(x_0I_n+X_0,x_1I_n+X_1)$ with scalars $x_0,x_1\in \K$ and traceless matrices
$X_0,X_1\in \Sl(n,\K)$.

\begin{thm}
The Lie algebra 
$\Gl(n,\K)\ltimes\Gl(n,\K)_{ab}$
with the following associative commutative product
\begin{align*}
(a_0I_n+A_0,\;& a_1I_n+A_1)\cdot (b_0I_n+B_0,\;b_1I_n+B_1) \\[3mm]
=
\Big(& a_0b_0 I_n   \;\;+\;\; a_0B_0 + b_0A_0\;,\;\\
&\big(a_0b_1+a_1b_0+\text{tr}(A_0B_1+A_1B_0)-\tfrac2{n}\text{tr}(A_0B_0)\big)I_n \;\;+\;\;
 a_0B_1 + b_0A_1+ A_0B_0+B_0A_0
 \Big)
\end{align*}
(here $A_0,A_1,B_0,B_1\in \Sl(n,\K)$)
is a Poisson algebra. 
For $n=3$,  it is isomorphic to $H_{E}^{\bullet,\bullet}(\h)$.
The $\Gl(n,\K)$-tables of the these products, with respect to 
the decomposition \eqref{eq.decomp_semidirect} and the partial labeling of Example \ref{ex.G=GL_k}, are
\noindent \begin{center}
{\small{}}%
\begin{tabular}{|c||c|c|c|c|}
\hline 
{\small{}$\cdot$} & {\small{}$(I_{n})_0$} & {\small{}$\Sl(n)_{0}$} & {\small{}$\Sl(n)_{ab}$} & {\small{}$(I_{n})_{ab}$}\tabularnewline
\hline 
\hline 
{\small{}$(I_{n})_0$} & {\small{}$(I_{n})_0$} & {\small{}$\Sl(n)_{0}$} & {\small{}$\Sl(n)_{ab}$} & {\small{}$(I_{n})_{ab}$}\tabularnewline
\hline 
{\small{}$\Sl(n)_{0}$} & {\small{}$\Sl(n)_{0}$} & {\small{}$\Sl(n)_{ab}^{(+)}$} & {\small{}$(I_{n})_{ab}$} & {\small{}$0$}\tabularnewline
\hline 
{\small{}$\Sl(n)_{ab}$} & {\small{}$\Sl(n)_{ab}$} & {\small{}$(I_{n})_{ab}$} & {\small{}$0$} & {\small{}$0$}\tabularnewline
\hline 
{\small{}$(I_{n})_{ab}$} & {\small{}$(I_{n})_{ab}$} & {\small{}$0$} & {\small{}$0$} & {\small{}$0$}\tabularnewline
\hline 
\end{tabular}
\end{center}
\begin{center}
\begin{tabular}{|c||c|c|c|c|}
\hline 
{\small{}$\left\{ -,-\right\} $} & {\small{}$(I_{n})_0$} & {\small{}$\Sl(n)_{0}$} & {\small{}$\Sl(n)_{ab}$} & {\small{}$(I_{n})_{ab}$}\tabularnewline
\hline 
\hline 
{\small{}$(I_{n})_0$} & {\small{}$0$} & {\small{}$0$} & {\small{}$0$} & {\small{}$0$}\tabularnewline
\hline 
{\small{}$\Sl(n)_{0}$} & {\small{}$0$} & {\small{}$\Sl(n)_{0}^{(-)}$} & {\small{}$\Sl(n)_{ab}^{(-)}$} & {\small{}$0$}\tabularnewline
\hline 
{\small{}$\Sl(n)_{ab}$} & {\small{}$0$} & {\small{}$\Sl(n)_{ab}^{(-)}$} & {\small{}$0$} & {\small{}$0$}\tabularnewline
\hline 
{\small{}$(I_{n})_{ab}$} & {\small{}$0$} & {\small{}$0$} & {\small{}$0$} & {\small{}$0$}\tabularnewline
\hline 
\end{tabular}{\small\par}
\par\end{center}
The $(-)$ and the $(+)$ indicate that the multiplications are $m_1^{Ad,Ad,Ad}$ and $m_2^{Ad,Ad,Ad}$ respectively (see Example \ref{ex.G=GL_k}).
\end{thm}

\begin{proof}
First, we must show that the given product and bracket 
define a Poisson algebra. 
These are lengthy but easy computations: the associativity 
of the product follows from that of matrix multiplication, 
and the Leibniz's law follows from the `associativity rule' 
$\text{tr}([A,B]C])=\text{tr}(A[B,C])$ 
and from the Leibniz's law $[A,BC]=B[A,C]+[A,B]C$.
The $\Gl(n,\K)$-tables are straightforward from the definition of 
the associative an Lie structure of 
$\Gl(n,\K)\ltimes\Gl(n,\K)_{ab}$ (see Example \ref{ex.G=GL_k}) 
and that this Poisson algebra is, for $n=3$ isomorphic to 
$H_{E}^{\bullet,\bullet}(\h)$ follows from 
Corollary \ref{Thm.Q3-HE(h1)}.
\end{proof}

\bibliographystyle{amsalpha}
\bibliography{bibliografia}

\end{document}